\documentclass{article}%
\usepackage{amssymb}
\usepackage{amsmath}
\usepackage{amsfonts}
\usepackage{graphicx}%
\providecommand{\U}[1]{\protect\rule{.1in}{.1in}}
\newtheorem{theorem}{Theorem}

\newtheorem{corollary}[theorem]{Corollary}

\newtheorem{lemma}[theorem]{Lemma}

\begin{document}

\title{Multiple Testing via Relative Belief Ratios}
\author{Michael Evans and Jabed Tomal\\Department of Statistical Sciences\\University of Toronto}
\date{}
\maketitle

\noindent\textit{Abstract}: Some large scale inference problems are considered
based on using the relative belief ratio as a measure of statistical evidence.
This approach is applied to the multiple testing problem. A particular
application of this is concerned with assessing sparsity. The approach taken
to sparsity has several advantages as it is based on a measure of evidence and
does not require that the prior be restricted in any way.\medskip

\noindent\textit{Key words and phrases}: multiple testing, sparsity,
statistical evidence, relative belief ratios, priors, checking for prior-data
conflict, testing for sparsity, relative belief multiple testing algorithm.

\section{Introduction}

Suppose there is a statistical model $\{f_{\theta}:\theta\in\Theta\}$ and the
parameter of interest is given by $\Psi:\Theta\rightarrow\Psi$ (we don't
distinguish between the function and its range to save notation) where $\Psi$
is an open subset of $R^{k}.$ Our concern here is with assessing individual
hypotheses $H_{0i}=\{\theta:\Psi_{i}(\theta)=\psi_{0i}\},$ namely, $H_{0i}$ is
the hypothesis that the $i$-th coordinate of $\psi$ equals $\psi_{0i}.$ If it
is known that $\psi_{i}=\Psi_{i}(\theta)=\psi_{0i},$ then the effective
dimension of the parameter of interest is $k-1$ which results in a
simplification of the model. For example, the $\psi_{i}$ may be linear
regression coefficients and $\psi_{0i}=0$ is equivalent to saying that an
individual term in the regression is dropped and so simplifies the relationship.

Considering all these hypotheses separately is the multiple testing problem
and the concern is to ensure that, while controlling the individual error
rate, the overall error rate does not become too large especially when $k$ is
large. An error means either the acceptance of $H_{0i}$ when it is false (a
false negative) or the rejection of $H_{0i}$ when it is true (a false
positive). This problem is considered here from a Bayesian perspective. A
first approach is to make an inference about the number of $H_{0i}$ that are
true (or false) and then use this to control the number of $H_{0i}$ that are
accepted (or rejected). Problems arise when $k$ is large, however, due to an
issue that arises with typical choices made for the prior. This paper
considers these problems and proposes solutions.

The problem of making inferences about sparsity is clearly related. For
suppose that there is a belief that many of the hypotheses $H_{0i}$ are indeed
true but there is little prior knowledge about which of the $H_{0i}$ are true.
This is effectively the sparsity problem. In this case it is not clear how to
choose a prior that reflects this belief. A common approach in the regression
problem is to use a prior that, together with a particular estimation
procedure, forces many of the $\psi_{i}$ to take the value $\psi_{0i}.$ A
difficulty with such an approach is that it is possible that such an
assignment is an artifact of the prior and the estimation procedure and not
primarily data-based. For example, the use of a Laplace prior together with
MAP\ estimation is known to accomplish this in certain problems. It would be
preferable, however, to have a procedure that was not dependent on a specific
form for the prior but rather worked with any prior and was based on the
statistical evidence for such an assignment contained in the data. The
methodology for multiple testing developed here accomplishes this goal.

The prior $\pi$ on $\theta$ clearly plays a key role in our developments.
Suppose that $\pi$ is decomposed as $\pi(\theta)=\pi(\theta\,|\,\psi)\pi
_{\Psi}(\psi)$ where $\pi(\theta\,|\,\psi)$ is the conditional prior on
$\theta$ given $\Psi(\theta)=\psi$ and $\pi_{\Psi}$ is the marginal prior on
$\psi.$ If $\pi_{\Psi}$ is taken to be too concentrated about $\psi_{0},$ then
there is the risk of encountering prior-data conflict whenever some of the
$H_{0i}$ are false. Prior-data conflict arises whenever the true value of the
parameter lies in a region of low prior probability, such as the tails of the
prior, and can lead to erroneous inferences. Consistent methods for addressing
prior-data conflict are developed in Evans and Moshonov (2006), Evans and Jang
(2011a) and a methodology for modifying the prior, when a conflict is
detected, is presented in Evans and Jang (2011b). Generally, a check for
prior-data conflict is made to ensure that the prior has been chosen
reasonably, just as model checking is performed to ensure that the model is
reasonable, given the observed data. The avoidance of prior-data conflict
plays a key role in our developments.

The assessment of the evidence for the truth of $H_{0i}$ is based on a measure
of the evidence that $H_{0i}$ is true as given by the relative belief ratio
for $\psi_{i}.$ Relative belief ratios are similar to Bayes factors and, when
compared to p-values, can be considered as more appropriate measures of
statistical evidence. For example, evidence can be obtained either for or
against a hypothesis. Moreover, there is a clear assessment of the strength of
this evidence so that when weak evidence is obtained, either for or against a
hypothesis, this does not entail acceptance or rejection, respectively.
Relative belief ratios and the theory of inference based on these quantities,
is discussed in Evans (2015) and a general outline is given in Section 2.
Section 3 considers applying relative belief to the multiple testing problem.
In Section 4 some practical applications are presented with special emphasis
on regression problems including the situation where the number of predictors
exceeds the number of observations.

There have been a number of priors proposed for the sparsity problem in the
literature. A common choice is the spike-and-slab prior due to George and
McCulloch (1993). Other priors are available that are known to produce
sparsity at least in connection with MAP estimation. For example, the Laplace
prior has a close connection with sparsity via the LASSO, as discussed in
Tibshirani (1996), Park and Casella (2008) and Hastie, Tibshirani and
Wainwright (2015). The horseshoe prior of Carvalho, Polson and Scott (2009)
will also produce sparsity through estimation. Both the Laplace and the
horseshoe prior can be checked for prior-data conflict but the check is not
explicitly for sparsity, only that the location and scaling of the prior is
appropriate. This reinforces a comment in Carvalho, Polson and Scott (2009)
that the use of the spike-and-slab prior represents a kind of "gold standard"
for sparsity. The problem with the spike-and-slab prior, and not possessed to
the same degree by the Laplace or horseshoe prior, is the difficulty of the
computations when implementing the posterior analysis. Basically there are
$2^{k}$ mixture components to the posterior and this becomes computationally
intractable as $k$ rises. Various approaches can be taken to try and deal with
this issue such as those developed in Rockova and George (2014).\ Any of these
priors can be used with the approach taken here but our methodology does not
require that a prior induces sparsity in any sense. Indeed sparsity is induced
only when the evidence, as measured by how the data changes beliefs, points to
this. This is a strong point of the approach as sparsity arises only through
an assessment of the evidence.

\section{Inferences Based on Relative Belief Ratios}

Suppose now that, in addition to the statistical model $\{f_{\theta}:\theta
\in\Theta\}$ there is a prior $\pi$ on $\Theta.$ After observing data $x,$ the
posterior distribution of $\theta$ is then given by the density $\pi
(\theta\,|\,x)=\pi(\theta)f_{\theta}(x)/m(x)$ where $m(x)=\int_{\Theta}%
\pi(\theta)f_{\theta}(x)\,d\theta$ is the prior predictive density of $x.$ For
a parameter of interest $\Psi(\theta)=\psi$ denote the prior and posterior
densities of $\psi$ by $\pi_{\Psi}$ and $\pi_{\Psi}(\cdot\,|\,x),$
respectively. The \textit{relative belief ratio} for a value $\psi$ is then
defined by $RB_{\Psi}(\psi\,|\,x)=\lim_{\delta\rightarrow0}\Pi_{\Psi
}(N_{\delta}(\psi\,)|\,x)/$\newline$\Pi_{\Psi}(N_{\delta}(\psi\,))$ where
$N_{\delta}(\psi\,)$ is a sequence of neighborhoods of $\psi$ converging
(nicely) to $\{\psi\}$ as $\delta\rightarrow0.$ Quite generally
\begin{equation}
RB_{\Psi}(\psi\,|\,x)=\pi_{\Psi}(\psi\,|\,x)/\pi_{\Psi}(\psi), \label{relbel}%
\end{equation}
the ratio of the posterior density to the prior density at $\psi.$ So
$RB_{\Psi}(\psi\,|\,x)$ is measuring how beliefs have changed concerning
$\psi$ being the true value from \textit{a priori} to \textit{a posteriori}.
When $RB_{\Psi}(\psi\,|\,x)>1$ the data have lead to an increase in the
probability that $\psi$ is correct and so there is evidence in favor of
$\psi,$ when $RB_{\Psi}(\psi\,|\,x)<1$ the data have lead to a decrease in the
probability that $\psi$ is correct and so there is evidence against $\psi,$
and when $RB_{\Psi}(\psi\,|\,x)=1$ there is no evidence either way. While
there are numerous quantities that measure evidence through change in belief,
the relative belief ratio is perhaps the simplest and it possesses many nice
properties as discussed in Evans (2015). For example, $RB_{\Psi}(\psi\,|\,x)$
is invariant under smooth changes of variable and also invariant to the choice
of the support measure for the densities. As such all relative belief
inferences possess this invariance which is not the case for many Bayesian
inferences such as using a posterior mode (MAP)\ or expectation for estimation.

The value $RB_{\Psi}(\psi_{0}\,|\,x)$ measures the evidence for the hypothesis
$H_{0}=\{\theta:\Psi(\theta)=\psi_{0}\}.$ It is also necessary to calibrate
whether this is strong or weak evidence for or against $H_{0}.$ Certainly the
bigger $RB_{\Psi}(\psi_{0}\,|\,x)$ is than 1, the more evidence there is in
favor of $\psi_{0}$ while the smaller $RB_{\Psi}(\psi_{0}\,|\,x)$ is than 1,
the more evidence there is against $\psi_{0}.$ But what exactly does a value
of $RB_{\Psi}(\psi_{0}\,|\,x)=20$ mean? It would appear to be strong evidence
in favor of $\psi_{0}$ because beliefs have increased by a factor of 20 after
seeing the data. But what if other values of $\psi$ had even larger increases?
A\ calibration of (\ref{relbel}) is then given by the \textit{strength}%
\begin{equation}
\Pi_{\Psi}(RB_{\Psi}(\psi\,|\,x)\leq RB_{\Psi}(\psi_{0}\,|\,x)\,|\,x),
\label{strength}%
\end{equation}
namely, the posterior probability that the true value of $\psi$ has a relative
belief ratio no greater than that of the hypothesized value $\psi_{0}.$ While
(\ref{strength}) may look like a p-value, it has a very different
interpretation. For when $RB_{\Psi}(\psi_{0}\,|\,x)<1,$ so there is evidence
against $\psi_{0},$ then a small value for (\ref{strength}) indicates a large
posterior probability that the true value has a relative belief ratio greater
than $RB_{\Psi}(\psi_{0}\,|\,x)$ and so there is strong evidence against
$\psi_{0}$ while only weak evidence against if (\ref{strength}) is big. If
$RB_{\Psi}(\psi_{0}\,|\,x)>1,$ so there is evidence in favor of $\psi_{0},$
then a large value for (\ref{strength}) indicates a small posterior
probability that the true value has a relative belief ratio greater than
$RB_{\Psi}(\psi_{0}\,|\,x))$ and so there is strong evidence in favor of
$\psi_{0},$ while a small value of (\ref{strength}) only indicates weak
evidence in favor of $\psi_{0}.$ Notice that in $\{\psi:RB_{\Psi}%
(\psi\,|\,x)\leq RB_{\Psi}(\psi_{0}\,|\,x)\},$ the `best' estimate of $\psi$
is given by $\psi_{0}$ because the evidence for this value is the largest in
the set.

Various results have been established in Baskurt and Evans (2013) and Evans
(2015) supporting both (\ref{relbel}), as the measure of the evidence, and
(\ref{strength}) as a measure of the strength of that evidence. For example,
the following simple inequalities, see Baskurt and Evans (2013), are useful in
assessing the strength as computing (\ref{strength}) can be avoided, namely,
$\Pi_{\Psi}(\{\psi_{0}\}\,|\,x)\leq\Pi_{\Psi}(RB_{\Psi}(\psi\,|\,x)\leq
RB_{\Psi}(\psi_{0}\,|\,x)\,|\,x)\leq RB_{\Psi}(\psi_{0}\,|\,x).$ So if
$RB_{\Psi}(\psi_{0}\,|\,x)>1$ and $\Pi_{\Psi}(\{\psi_{0}\}\,|\,x)$ is large,
there is strong evidence in favor of $\psi_{0}$ while, if $RB_{\Psi}(\psi
_{0}\,|\,x)<1$ is very small, then there is immediately strong evidence
against $\psi_{0}.$

The use of $RB_{\Psi}(\psi_{0}\,|\,x)$ to assess the hypothesis $H_{0}$ also
possesses optimality properties. For example, let $A$ be a subset of the
sample space such that whenever $x\in A,$ the hypothesis is accepted. If
$M(A\,|\,H_{0})$ denotes the prior predictive probability of $A$ given that
$H_{0}$ is true, then $M(A\,|\,H_{0})$ is the prior probability of accepting
$H_{0}$ when it is true. The relative belief acceptance region is naturally
given by $A_{rb}(\psi_{0})=\{x:RB_{\Psi}(\psi_{0}\,|\,x)>1\}.$ Similarly, let
$R$ be a subset such that whenever $x\in R,$ the hypothesis is rejected and
let $M(R\,|\,H_{0}^{c})$ denote the prior predictive probability of $R$ given
that $H_{0}$ is false. The relative belief rejection region is then naturally
given by $R_{rb}(\psi_{0})=\{x:RB_{\Psi}(\psi_{0}\,|\,x)<1\}.$ The following
result is proved in Evans (2015) as Proposition 4.7.9.

\begin{theorem}
\label{thm1}(i) The acceptance region $A_{rb}(\psi_{0})$ minimizes $M(A)$
among all acceptance regions $A$ satisfying $M(A\,|\,H_{0})\geq M(A_{rb}%
(\psi_{0})\,|\,H_{0}).$ (ii) The rejection region $R_{rb}(\psi_{0})$ maximizes
$M(R)$ among all rejection regions $R$ satisfying $M(R\,|\,H_{0})\leq
M(R_{rb}(\psi_{0})\,|\,H_{0}).$
\end{theorem}

\noindent To see the meaning of this result note that
\begin{align*}
&  M(A)=E_{\Pi_{\Psi}}(M(A\,|\,\psi))=E_{\Pi_{\Psi}}(I_{H_{0}^{c}}%
(\psi)M(A\,|\,\psi))+\Pi_{\Psi}(H_{0})M(A\,|\,H_{0})\\
&  \geq M(A_{rb}(\psi_{0}))=E_{\Pi_{\Psi}}(I_{H_{0}^{c}}(\psi)M(A_{rb}%
(\psi_{0})\,|\,\psi))+\Pi_{\Psi}(H_{0})M(A_{rb}(\psi_{0})\,|\,H_{0}).
\end{align*}
Therefore, if $\Pi_{\Psi}(H_{0})=0,$ then $A_{rb}(\psi_{0})$ minimizes
$E_{\Pi_{\Psi}}(I_{\{\psi_{0}\}^{c}}(\psi)M(A\,|\,\psi)),$ the prior
probability that $H_{0}$ is accepted given that it is false, among all
acceptance regions $A$ satisfying $M(A\,|\,\psi_{0})\geq M(A_{rb}(\psi
_{0})\,|\,\psi_{0})$ and when $\Pi_{\Psi}(H_{0})>0,$ then $A_{rb}(\psi_{0})$
minimizes this probability among all acceptance regions $A$ satisfying
$M(A\,|\,H_{0})=M(A_{rb}(\psi_{0})\,|\,H_{0}).$ Also,%
\begin{align*}
&  M(R)=E_{\Pi_{\Psi}}(M(R\,|\,\psi))=E_{\Pi_{\Psi}}(I_{H_{0}^{c}}%
(\psi)M(R\,|\,\psi))+\Pi_{\Psi}(H_{0})M(R\,|\,H_{0})\\
&  \leq M(R_{rb}(\psi_{0}))=E_{\Pi_{\Psi}}(I_{H_{0}^{c}}(\psi)M(R_{rb}%
(\psi_{0})\,|\,\psi))+\Pi_{\Psi}(H_{0})M(R_{rb}(\psi_{0})\,|\,H_{0}).
\end{align*}
Therefore, if $\Pi_{\Psi}(\{\psi_{0}\})=0,$ then $R_{rb}(\psi_{0})$ maximizes
$E_{\Pi_{\Psi}}(I_{\{\psi_{0}\}^{c}}(\psi)M(R\,|\,\psi)),$ the prior
probability that $H_{0}$ is rejected given that it is false, among all
rejection regions $R$ satisfying $M(R\,|\,H_{0})\leq M(R_{rb}(\psi
_{0})\,|\,H_{0})$ and when $\Pi_{\Psi}(H_{0}\})>0,$ then $R_{rb}(\psi_{0})$
maximizes this probability among all rejection regions $R$ satisfying
$M(R\,|\,\psi_{0})=M(R_{rb}(\psi_{0})\,|\,\psi_{0}).$

Note that it does not make sense to accept or reject $H_{0}$ when $RB_{\Psi
}(\psi_{0}\,|\,x)=1.$ Also, under i.i.d. sampling, $M(A_{rb}(\psi
_{0})\,|\,H_{0})\rightarrow1$ and $M(R_{rb}(\psi_{0})\,|\,H_{0})\rightarrow0$
as sample size increases, so these quantities can be controlled by design.
When $\Pi_{\Psi}(\{\psi_{0}\})=0,$ then \ $M(A_{rb}(\psi_{0})\,|\,H_{0}%
)=1-M(R_{rb}(\psi_{0})\,|\,H_{0})$ and so controlling $M(A_{rb}(\psi
_{0})\,|\,H_{0})$ is controlling the "size" of the test. In general,
$E_{\Pi_{\Psi}}(I_{\{\psi_{0}\}^{c}}(\psi)M(R_{rb}(\psi_{0})\,|\,\psi))$ can
be thought of as the Bayesian power of the relative belief test. Note that it
is reasonable to set either the probability of a false positive or the
probability of a true negative as part of design and so the theorem is an
optimality result with practical import. It is easily seen that the proof
of\ Theorem 1 can be generalized to obtain optimality results for the
acceptance region $A_{rb,q}(\psi_{0})=\{x:RB_{\Psi}(\psi_{0}\,|\,x)>q\}$ and
for the rejection region $R_{rb,q}(\psi_{0})=\{x:RB_{\Psi}(\psi_{0}%
\,|\,x)<q\}.$ The following inequality is useful in\ Section 3 in controlling
error rates.

\begin{theorem}
\label{thm2}$M\left(  \left.  R_{rb,q}(\psi_{0})\,\right\vert \,\psi
_{0}\right)  \leq q.$
\end{theorem}

\noindent Proof: By the Savage-Dickey result, see Proposition 4.2.7 in Evans
(2015), $RB_{\Psi}(\psi_{0}\,|\,x)=m(x\,|\,\psi_{0})/m(x).$ Now $E_{M(\cdot
\,|\,\psi_{0})}(m(x)/m(x\,|\,\psi_{0}))=1$ and so by Markov's inequality,
$M\left(  \left.  R_{rb,q}(\psi_{0})\,\,\right\vert \,\psi_{0}\right)
=M\left(  \left.  m(x)/m(x\,|\,\psi_{0})>1/q\,\right\vert \,\psi_{0}\right)
\leq q.$\smallskip

There is another issue associated with using $RB_{\Psi}(\psi_{0}\,|\,x)$ to
assess the evidence that $\psi_{0}$ is the true value. One of the key concerns
with Bayesian inference methods is that the choice of the prior can bias the
analysis in various ways. \ For example, in many problems Bayes factors and
relative belief ratios can be made arbitrarily large by choosing the prior to
be increasingly diffuse. This phenomenon is associated with the
Jeffreys-Lindley paradox and clearly indicates that it is possible to bias
results by the choice of the prior.

An approach to dealing with this bias is discussed in Baskurt and Evans
(2013). For, given a measure of evidence, it is possible to measure \textit{a
priori} whether or not the chosen prior induces bias either for or against
$\psi_{0}.$ The bias against $\psi_{0}$ is given by
\begin{equation}
M(\left.  RB_{\Psi}(\psi_{0}\,|\,x)\leq1\,\right\vert \,\psi_{0}%
)=1-M(A_{rb}(\psi_{0})\,|\,\psi_{0}) \label{bias1}%
\end{equation}
as this is the prior probability that evidence will not be obtained in favor
of $\psi_{0}$ when $\psi_{0}$ is true. If (\ref{bias1}) is large, subsequently
reporting, after seeing the data, that there is evidence against $\psi_{0}$ is
not convincing.

To measure the bias in favor of $\psi_{0}$ choose values $\psi_{0}^{\prime
}\neq\psi_{0}$ such that the difference between $\psi_{0}$ and $\psi
_{0}^{\prime}$ represents the smallest difference of practical importance.
Then compute
\begin{equation}
M_{T}\left(  \left.  RB_{\Psi}(\psi_{0}\,|\,x)\geq1\,\right\vert \,\psi
_{0}^{\prime}\right)  =1-M(R_{rb}(\psi_{0})\,|\,\psi_{0}^{\prime})
\label{bias2}%
\end{equation}
as this is the prior probability that we will not obtain evidence against
$\psi_{0}$ when $\psi_{0}$ is false. Note that (\ref{bias2}) tends to decrease
as $\psi_{0}^{\prime}$ moves further away from $\psi_{0}.$ When (\ref{bias2})
is large, subsequently reporting, after seeing the data, that there is
evidence in favor of $\psi_{0}$, is not convincing. For a fixed prior, both
(\ref{bias1}) and (\ref{bias2}) decrease with sample size and so, in design
situations, they can be used to set sample size and so control bias. Notice
that $M(A_{rb}(\psi_{0})\,|\,\psi_{0})$ can be considered as the sensitivity
and $M(R_{rb}(\psi_{0})\,|\,\psi_{0}^{\prime})$ as the specificity of the
relative belief hypothesis assessment. These issues are further discussed in
Evans (2015).

As $RB_{\Psi}(\psi\,|\,x)$ measures the evidence that $\psi$ is the true
value, it can also be used to estimate $\psi.$ For example, the best estimate
of $\psi$ is clearly the value for which the evidence is greatest, namely,
$\psi(x)=\arg\sup RB_{\Psi}(\psi\,|\,x).$ Associated with this is a $\gamma
$-credible region $C_{\Psi,\gamma}(x)=\{\psi:RB_{\Psi}(\psi\,|\,x)\geq
c_{\Psi,\gamma}(x)\}$ where $c_{\Psi,\gamma}(x)=\inf\{k:\Pi_{\Psi}(RB_{\Psi
}(\psi\,|\,x)>k\,|\,x)\leq\gamma\}.$ Notice that $\psi(x)\in C_{\Psi,\gamma
}(x)$ for every $\gamma\in\lbrack0,1]$ and so, for selected $\gamma,$ we can
take the "size" of $C_{\Psi,\gamma}(x)$ as a measure of the accuracy of the
estimate $\psi(x).$ The interpretation of $RB_{\Psi}(\psi\,|\,x)$ as the
evidence for $\psi\ $\ forces the sets $C_{\Psi,\gamma}(x)$ to be the credible
regions. For if $\psi_{1}$ is in such a region and $RB_{\Psi}(\psi
_{2}\,|\,x)\geq RB_{\Psi}(\psi_{1}\,|\,x),$ then $\psi_{2}$ must be in the
region as well as there is at least as much evidence for $\psi_{2}$ as for
$\psi_{1}.$ A variety of optimality results have been established for
$\psi(x)$ and $C_{\Psi,\gamma}(x),$ see Evans (2015).

The view is taken here that any time continuous probability is used, then this
is an approximation to a finite, discrete context. For example, if $\psi$ is a
mean and the response measurements are to the nearest centimeter, then of
course the true value of $\psi$ cannot be known to an accuracy greater than
1/2 of a centimeter, no matter how large a sample we take. Furthermore, there
are implicit bounds associated with any measurement process. As such the
restriction is made here to discretized parameters that take only finitely
many values. So when $\psi$ is a continuous, real-valued parameter, it is
discretized to the intervals $\ldots,(\psi_{0}-3\delta,\psi_{0}-\delta
],(\psi_{0}-\delta,\psi_{0}+\delta],(\psi_{0}+\delta,\psi_{0}+3\delta],\ldots$
for some choice of $\delta>0,$ and there are only finitely may such intervals
covering the range of possible values. It is of course possible to allow the
intervals to vary in length as well. With this discretization, then we can
take $H_{0}=$ $(\psi_{0}-\delta,\psi_{0}+\delta].$

\section{Inferences for Multiple Tests}

Consider now the multiple testing problem discussed in Section 1. Let $\xi
=\Xi(\theta)=\frac{1}{k}\sum_{i=1}^{k}I_{H_{0i}}(\Psi_{i}(\theta))$ be the
proportion of the hypotheses $H_{0i}$ that are true and suppose that $\Psi
_{i}$ is finite for each $i,$ perhaps arising via a discretization as
discussed in Section 2. Note that the discreteness is essential to
realistically determine what proportion of the hypotheses are correct,
otherwise, under a continuous prior on $\Psi,$ the prior distribution of
$\Xi(\theta)$ is degenerate at 0. In an application it is desirable to make
inference about the true value of $\xi\in\Xi=\{0,1/k,2/k,\ldots,1\}.$ This is
based on the relative belief ratio $RB_{\Xi}(\xi\,|\,\,x)=\Pi(\Xi(\theta
)=\xi\,|\,x)/\Pi(\Xi(\theta)=\xi).$ The appropriate estimate of $\xi$ is then
the relative belief estimate of $\Xi,$ namely, $\xi(x)=\arg\sup_{\xi}RB_{\Xi
}(\xi\,|\,x).$ The accuracy of this estimate is assessed using the size of
$C_{\Xi,\gamma}(x)$ for some choice of $\gamma\in\lbrack0,1].$ Also,
hypotheses such as $H_{0}=\{\theta:\Xi(\theta)\in\lbrack\xi_{0},\xi_{1}]\},$
namely, the proportion true is at least $\xi_{0}$ and no greater than $\xi
_{1},$ can be assessed using the relative belief ratio $RB(H_{0}%
\,|\,x)=\Pi(\xi_{0}\leq\Xi(\theta)\leq\xi_{1}\,|\,x)/\Pi(\xi_{0}\leq\Xi
(\theta)\leq\xi_{1})$ which equals $RB_{\Xi}(\xi_{0}\,|\,x)$ when $\xi_{0}%
=\xi_{1}.$ The strength of this evidence can be assessed as previously discussed.

The estimate $\xi(x)$ can be used to control how many hypotheses are
potentially accepted. For this select $k\xi(x)$ of the $H_{0i}$ as being true
from among those for which $RB_{\Psi_{i}}(\psi_{0i}\,|\,x)>1.$ Note that it
does not make sense to accept $H_{0i}$ as being true when $RB_{\Psi_{i}}%
(\psi_{0i}\,|\,x)<1$ as there is evidence against this hypothesis. So, if
there are fewer than $k\xi(x)$ satisfying $RB_{\Psi_{i}}(\psi_{0i}%
\,|\,x)>1,$\ then fewer than this number should be accepted. If there are more
than $k\xi(x)$ of the relative belief ratios satisfying $RB_{\Psi_{i}}%
(\psi_{0i}\,|\,x)>$\ $1,$ then some method will have to be used to select the
$k\xi(x)$ which are potentially accepted. It is clear, however, that the
logical way to do this is to order the $H_{0i},$ for which $RB_{\Psi_{i}}%
(\psi_{0i}\,|\,x)>$\ $1,$ based on their strengths $\Pi_{\Psi}(RB_{\Psi_{i}%
}(\psi_{0i}\,|\,x)\leq RB_{\Psi_{i}}(\psi_{0i}\,|\,x)\,|\,x),$ from largest to
smallest, and accept at most the $k\xi(x)$ for which the evidence is strongest.

Note too that, if some of these strengths are indeed weak, there is no need to
necessarily accept these hypotheses. The ultimate decision as to whether or
not to accept a hypothesis is application dependent and is not statistical in
nature. The role of statistics is to supply a clear statement of the evidence
and its strength, while other criteria come into play when decisions are made.
In any case, it is seen that inference about $\xi$ is being used to help
control the number of hypotheses accepted.

If, as is more common, control is desired of the number of false positives,
then the relevant parameter of interest is $\upsilon=\Upsilon(\theta
)=1-\Xi(\theta),$ the proportion of false hypotheses. Note that $\Pi
(\Upsilon(\theta)=\upsilon)=\Pi(\Xi(\theta)=1-\upsilon)$ and so the relative
belief estimate of $\upsilon$ satisfies $\upsilon(x)=1-\xi(x).$ Following the
same procedure, the $H_{0i}$ with $RB_{\Psi_{i}}(\psi_{0i}\,|\,x)<1$ are then
ranked via their strengths and at most $k\upsilon(x)$ are rejected.

The consistency of the procedure just described, for correctly identifying the
$H_{0i}$ that are true and those that are false, follows from results proved
in Section 4.7.1 of Evans (2015) under i.i.d. sampling. In other words, as the
amount of data increases, $\xi(x)$ converges to the proportion of $H_{0i}$
that are true, each $RB(\psi_{0i}\,|\,\,x)$ converges to the largest possible
value (always bigger than 1) when $H_{0i}$ is true and converges to 0 when
$H_{0i}$ is false, and the evidence in favor or against converges to the
strongest possible, depending on whether the hypothesis in question is true or false.

We refer to this procedure as the \textit{multiple testing algorithm.}
Consider first a somewhat artificial example where many of the computations
are easy but which fully demonstrates the relevant characteristics of the
methodology.\smallskip

\noindent\textbf{Example 1.} \textit{Location normal.}

Suppose that there are $k$ independent samples $x_{ij}$ for $1\leq i\leq
k,1\leq j\leq n$ where the $i$-th sample is from a $N(\mu_{i},\sigma^{2})$
distribution with $\mu_{i}$ unknown and $\sigma^{2}$ known. It also assumed
that prior knowledge about all the unknown $\mu_{i}$ is reflected in the
statement that the $\mu_{1},\ldots,\mu_{k}$ are i.i.d. from a $N(\mu
_{0},\lambda_{0}^{2}\sigma^{2})$ distribution. It is easy to modify our
development to allow the sample sizes to vary and to use a general
multivariate normal prior, while the case where $\sigma^{2}$ is unknown is
considered in Section 4. This context is relevant to the analysis of
microarray data.

The value of $(\mu_{0},\lambda_{0}^{2})$ is determined via elicitation. For
this it is supposed that it is known with virtual certainty that each $\mu
_{i}\in(m_{l},m_{u})$ for specified values $m_{l}\leq m_{u}.$ Here virtual
certainty is interpreted to mean that the prior probability of this interval
is at least $0.99$ and other choices could be made. It is also supposed that
$\mu_{0}=$ $(m_{l}+m_{u})/2.$ This implies that $\lambda_{0}=(m_{u}%
-m_{l})/(2\sigma\Phi^{-1}((1+0.99)/2)).$ Following Evans and Jang (2011b),
increasing the value of $\lambda_{0}$ implies a more weakly informative prior
in this context and, as such, decreases the possibility of prior-data conflict.

The posterior distributions of the $\mu_{i}$ are then independent with
$\mu_{i}\,|\,x\sim N(\mu_{i}(x),(n\lambda_{0}^{2}+1)^{-1}\lambda_{0}^{2}%
\sigma^{2})$ where $\mu_{i}(x)=(n+1/\lambda_{0}^{2})^{-1}(n\bar{x}_{i}+\mu
_{0}/\lambda_{0}^{2}).$ Given that the measurements are taken to finite
accuracy, it is not necessarily realistic to test $\mu_{i}=\mu_{0}$. As such,
a value $\delta>0$ is specified so that $H_{0i}=(\mu_{0}-\delta/2,\mu
_{0}+\delta/2]$ in a discretization of the parameter space into a finite
number of \ intervals, each of length $\delta,$ as well as two tail intervals.
Then for $T\in%
\mathbb{N}
$ there are $2T+1$ intervals of the form $(\mu_{0}+(t-1/2)\delta,\mu
_{0}+(t+1/2)\delta],$ for $t\in\{-T,-T+1,\ldots,T\}$ that span $(m_{l}%
,m_{u}),$ together with two additional tail intervals $(-\infty,\mu
_{0}-(T+1/2)\delta]$ and $(\mu_{0}+(T+1/2)\delta,\infty)$ to cover the full
range. The relative belief ratio for the $t$-th interval for $\mu_{i}$\ is
then given by
\begin{align}
&  RB_{_{i}}((\mu_{0}+(t-1/2)\delta,\mu_{0}+(t+1/2)\delta]\,|\,x)\nonumber\\
&  =\frac{\left\{
\begin{array}
[c]{c}%
\Phi((n\lambda_{0}^{2}+1)^{1/2}(\mu_{0}+(t+1/2)\delta-\mu_{i}(x))/\lambda
_{0}\sigma)-\\
\Phi((n\lambda_{0}^{2}+1)^{1/2}(\mu_{0}+(t-1/2)\delta-\mu_{i}(x))/\lambda
_{0}\sigma)
\end{array}
\right\}  }{\Phi((t+1/2)\delta/\lambda_{0}\sigma)-\Phi((t-1/2)\delta
/\lambda_{0}\sigma)} \label{relbel2}%
\end{align}
with a similar formula for the tail intervals. When $\delta$ is small this is
approximated by the ratio of the posterior to prior densities of $\mu_{i}$
evaluated at $\mu_{0}+t\delta.$ Then $RB(H_{0i}\,|\,x)=RB_{_{i}}((\mu
_{0}-\delta/2,\mu_{0}+\delta/2]\,|\,x)$ gives the evidence for or against
$H_{0i}$ and the strength of this evidence is computed using the discretized
posterior distribution. Notice that (\ref{relbel2}) converges to $\infty$ as
$\lambda_{0}\rightarrow\infty$ and this is characteristic of other measures of
evidence such as Bayes factors. As discussed in Evans (2015), this is one of
the reasons why calibrating (\ref{relbel}) via (\ref{strength}) is necessary.

As previously noted, it is important to take account of the bias in the prior.
To simplify matters the continuous approximation is used here as this makes
little difference for the discussion of bias concerning inference about
$\mu_{i}$ (see Tables \ref{biastab1} and 3)$.$ The bias against $\mu_{i}%
=\mu_{0}$ equals%
\begin{equation}
M(RB_{i}(\mu_{0}\,|\,x)\leq1\,|\,\mu_{0})=2(1-\Phi(a_{n}(1))) \label{bias1ex1}%
\end{equation}
where
\[
a_{n}(q)=\left\{
\begin{array}
[c]{cc}%
\{(1+1/n\lambda_{0}^{2})\log((n\lambda_{0}^{2}+1)/q^{2})\}^{1/2}, & q^{2}\leq
n\lambda_{0}^{2}+1\\
0, & q^{2}>n\lambda_{0}^{2}+1.
\end{array}
\right.
\]
Note that (\ref{bias1ex1}) converges to $2(1-\Phi(1))=0.32$ as $\lambda
_{0}\rightarrow0$ and to 0 as $\lambda_{0}\rightarrow\infty$ and, for fixed
$\lambda_{0},$ converges to 0 as $n\rightarrow\infty.$ So there is never
strong bias against $\mu_{i}=\mu_{0}$ and this is as expected since the prior
is centered on $\mu_{0}.$ The bias in favor of $\mu_{i}=\mu_{0}$ is measured
by%
\begin{equation}
M(RB_{i}(\mu_{0}\,|\,x)\geq1\,|\,\mu_{0}\pm\delta/2)=\Phi(\sqrt{n}%
\delta/2\sigma+a_{n}(1))-\Phi(\sqrt{n}\delta/2\sigma-a_{n}(1)).
\label{bias2ex1}%
\end{equation}
As $\lambda_{0}\rightarrow\infty$ then (\ref{bias2ex1}) converges to 1 so
there is bias in favor of $\mu_{i}=\mu_{0}$ and this reflects what was
obtained for the limiting value of (\ref{relbel2}). Also this decreases with
increasing $\delta$ and goes to 0 as $n\rightarrow\infty.$ So indeed bias of
both types can be controlled by sample size. Perhaps the most important take
away from this discussion, however, is that by using a supposedly
noninformative prior with $\lambda_{0}$ large, bias in favor of the $H_{0i}$
is being induced.

Consider first a simulated data set $x$ when $k=10,n=5,\sigma=1,\delta
=1,\mu_{0}=0,(m_{l},m_{u})=(-5,5),$ so that $\lambda_{0}=10/2\Phi
^{-1}(0.995)=1.94$ and suppose $\mu_{1}=$ $\mu_{2}=\cdots=\mu_{7}=0$ with the
remaining $\mu_{i}=2.$ The relative belief ratio function $RB_{\Xi}%
(\cdot\,|\,\,x)$ is plotted in Figure 1. In this case the relative belief
estimate $\xi(x)=0.70$ is exactly correct. Table 1 gives the values of the
$RB_{i}(0\,|\,x)$ together with their strengths. It is clear that the multiple
testing algorithm leads to 0 false positives and 0 false negatives. So the
algorithm works perfectly on this data but of course it can't be expected to
do as well when the three nonzero means move closer to 0. Also, it is worth
noting that the strength of the evidence in favor of $\mu_{i}=0$ is very
strong for $i=1,2,3,5,6,7$ but only moderate when $i=4.$ The strength of the
evidence against $\mu_{i}=0$ is very strong for $i=8,9,10.$ Note that the
maximum possible value of $RB_{_{i}}((\mu_{0}-\delta/2,\mu_{0}+\delta
/2]\,|\,x)$ here is $(2\Phi(\delta/2\lambda_{0}\sigma)-1)^{-1}=4.92,$ so
indeed some of the relative belief ratios are relatively large.

Now consider basically the same context but where $k=1000$ and $\mu_{1}=$
$\cdots=\mu_{700}=0$ while the remaining 300 satisfy $\mu_{i}=2.$ The relative
belief ratio $RB_{\Xi}(\cdot\,|\,\,x)$ is plotted in Figure 2. In this case
$\xi(x)=0.47$ which is a serious underestimate. As such the multiple testing
algorithm will not record enough acceptances and so will fail.%

\begin{figure}
[ptb]
\begin{center}
\includegraphics[
height=2.2in,
width=2.4in
]%
{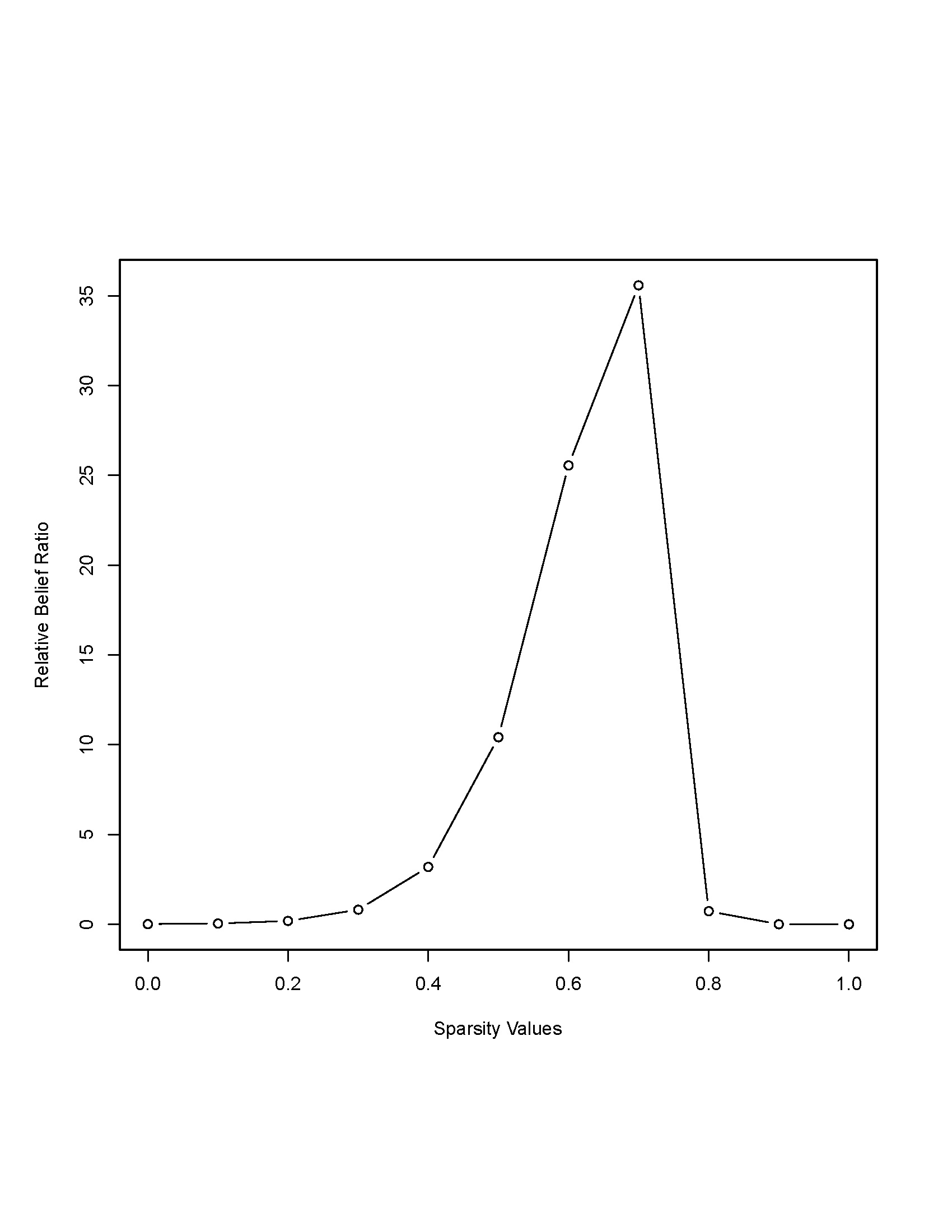}
\caption{A plot of the relative belief ratio of $\Xi$\ when $n=5,k=10$ and $7$
means equal 0 with the remaining means equal to $2$ in Example 1.}
\end{center}
\end{figure}

\begin{table}[tbp] \centering
\begin{tabular}
[c]{|l|c|c|c|c|c|}\hline
$i$ & $1$ & $2$ & $3$ & $4$ & $5$\\\hline
$\mu_{i}$ & $0$ & $0$ & $0$ & $0$ & $0$\\
$RB_{i}(0\,|\,x)$ & $3.27$ & $3.65$ & $2.98$ & $1.67$ & $3.57$\\
Strength & $1.00$ & $1.00$ & $1.00$ & $0.37$ & $1.00$\\\hline
$i$ & $6$ & $7$ & $8$ & $9$ & $10$\\\hline
$\mu_{i}$ & $0$ & $0$ & $2$ & $2$ & $2$\\
$RB_{i}(0\,|\,x)$ & $3.00$ & $3.43$ & $2.09\times10^{-4}$ & $3.99\times
10^{-4}$ & $8.80\times10^{-3}$\\
Strength & $1.00$ & $1.00$ & $4.25\times10^{-5}$ & $8.11\times10^{-5}$ &
$1.83\times10^{-3}$\\\hline
\end{tabular}
\caption{Relative belief ratios and strengths for the $\mu_i$ in Example 1 with
$k=10$.}\label{TableKey}%
\end{table}%

\begin{figure}
[ptb]
\begin{center}
\includegraphics[
height=2.2in,
width=2.4in
]%
{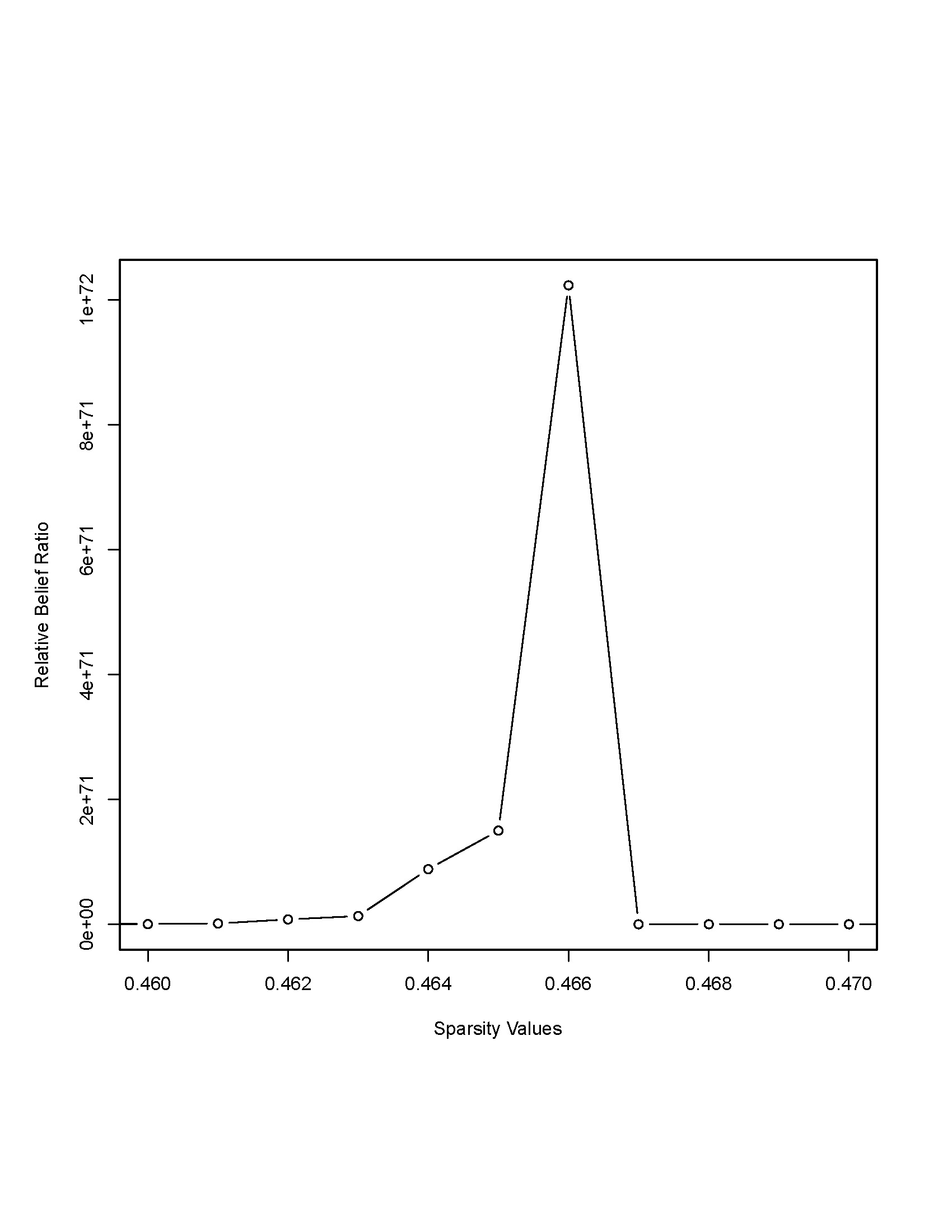}%
\caption{A plot of the relative
belief ratio of $\Xi$ when $n=5,k=1000$ and $700$ means equal $0$ with the
remaining means equal to $2$ in Example 1.}
\end{center}
\end{figure}

This problem arises due to the independence assumption on the $\mu_{i}.$ For
the prior distribution of $k\Xi(\theta)$ is binomial$(k,2\Phi(\delta
/2\lambda_{0}\sigma)-1)$ and the prior distribution of $k\Upsilon(\theta)$ is
binomial$(k,2(1-\Phi(\delta/2\lambda_{0}\sigma))).$ So the \textit{a priori}
expected proportion of true hypotheses is $2\Phi(\delta/2\lambda_{0}\sigma)-1$
and the expected proportion of false hypotheses is $2(1-\Phi(\delta
/2\lambda_{0}\sigma)).$ When $\delta/2\lambda_{0}\sigma$ is small, as when the
amount of sampling variability or the diffuseness of the prior are large, then
the prior on $\Xi$ suggests a belief in many false hypotheses. When $k$ is
small, relatively small amounts of data can override this to produce accurate
inferences about $\xi$ or $\upsilon$ but otherwise large amounts of data are
needed which may not be available. $\blacksquare\smallskip$

Given that accurate inference about $\xi$ and $\upsilon$ is often not
feasible, we focus instead on protecting against too many false positives and
false negatives. From the discussion in Example 1 it is seen that it is
possible to produce bias in favor of the $H_{0i}$ being true simply by using a
diffuse prior. If our primary goal is to guard against too many false
positives, then this biasing will work as we can make the prior conditional
probabilities of the events $RB(\psi_{0i}\,|\,\,x)<1,$ given that $H_{0i}$ is
true, as small as desirable by choosing such a prior. This could be seen as a
way of creating a higher bar for a positive result. The price we pay for this,
however, is too many false negatives. As such we consider another way of
"raising the bar". Given that $RB(\psi_{0i}\,|\,\,x)$ is measuring evidence, a
natural approach is to simply choose constants $0<q_{R}\leq1\leq q_{A}$ and
classify $H_{0i}$ as accepted when $RB(\psi_{0i}\,|\,\,x)>q_{A}$ and rejected
when $RB(\psi_{0i}\,|\,\,x)<q_{R}.$ The strengths can also be quoted to assess
the reliability of these inferences. Provided $q_{R}$ is greater than the
minimum possible value of $RB(\psi_{0i}\,|\,\,x),$ and this is typically 0,
and $q_{A}$ is chosen less than the maximum possible value of $RB(\psi
_{0i}\,|\,\,x)$, and this is 1 over the prior probability of $H_{0i},$ then
this procedure is consistent as the amount of data increases. In fact, the
related estimates of $\xi$ and $\upsilon$ are also consistent. The price paid
for this is that a hypothesis will not be classified whenever $q_{R}\leq
RB(\psi_{0i}\,|\,\,x)\leq q_{A}.$ Not classifying a hypothesis implies that
there is not enough evidence for this purpose and more data is required. This
approach is referred to as the \textit{relative belief multiple testing
algorithm}.

It remains to determine $q_{A}$ and $q_{R}$. Consider first protecting against
too many false positives. The \textit{a priori} conditional prior probability,
given that $H_{0i}$ is true, of finding evidence against $H_{0i}$ less than
$q_{R}$ is $M(RB_{i}(\psi_{0i}\,|\,X)<q_{R}\,|\,\psi_{0i})=M(m(X\,|\,\psi
_{0i})/m(x)<q_{R}\,|\,\psi_{0i})\leq q_{R}$ where the final inequality follows
from Theorem 2. Naturally, we want the probability of a false negative to be
small and so choosing $q_{R}$ small accomplishes this. The \textit{a priori
}probability that a randomly selected hypothesis produces a false positive is
\begin{equation}
\frac{1}{k}\sum_{i=1}^{k}M(RB_{i}(\psi_{0i}\,|\,X)<q_{R}\,|\,\psi_{0i})
\label{exfalspos}%
\end{equation}
which by Theorem \ref{thm2} is bounded above by $q_{R}$ and so converges to 0
as $q_{R}\rightarrow0.$ Also, for fixed $q_{R},$ (\ref{exfalspos}) converges
to 0 as the amount of data increases. More generally $q_{R}$ can be allowed to
depend on $i$ but when the $\psi_{i}$ are similar in nature this does not seem
necessary. Furthermore, it is not necessary to weight the hypotheses equally
so a randomly chosen hypothesis with unequal probabilities could be relevant
in \ certain circumstances. In any case, controlling the value of
(\ref{exfalspos}), whether by sample size or by the choice of $q_{R},$ is
clearly controlling for false positives. Suppose there is proportion $p_{FP}$
of false positives that is just tolerable in a problem. Then $q_{R}$ can be
chosen so that (\ref{exfalspos}) is less than or equal to $p_{FP}$ and note
that $q_{R}=p_{FP}$ satisfies this.

Similarly, if $\psi_{0i}^{\prime}\neq\psi_{0i}$ then $M(RB_{i}(\psi
_{0i}\,|\,X)>q_{A}\,|\,\psi_{0i}^{\prime})$ is the prior probability of
finding evidence for $H_{0i}$ when $\psi_{0i}^{\prime}$ is the true value. For
a given effect size $\delta$ of practical importance it is natural to take
$\psi_{0i}^{\prime}=\psi_{0i}\pm\delta/2.$ In typical applications this
probability becomes smaller the "further" $\psi_{0i}^{\prime}$ is from
$\psi_{0i}$ and so choosing $q_{A}$ to make this probability small will make
it small for all alternatives. Under these circumstances the \textit{a priori
}probability that a randomly selected hypothesis produces a false negative is
bounded above by%
\begin{equation}
\frac{1}{k}\sum_{i=1}^{k}M(RB_{i}(\psi_{0i}\,|\,X)>q_{A}\,|\,\psi_{0i}%
^{\prime}). \label{exfalsnegs}%
\end{equation}
As $q_{A}\rightarrow\infty,$ or as\thinspace the amount of data increases with
$q_{A}$ fixed, then (\ref{exfalsnegs}) converges to 0 so the number of false
negatives can be controlled. If there is proportion $p_{FN}$ of false
negatives that is just tolerable in a problem, then $q_{A}$ can be chosen so
that (\ref{exfalsnegs}) is less than or equal to $p_{FN}.$

The following optimality result holds for relative belief multiple testing.

\begin{corollary}
\label{cor3}(i) Among all procedures where the prior probability of accepting
$H_{0i}$, when it is true, is at least $M(RB_{i}(\psi_{0i}\,|\,X)>q_{A}%
\,|\,\psi_{0i})$ for $i=1,\ldots,k$, the relative belief multiple testing
algorithm minimizes the prior probability that a randomly chosen hypothesis is
accepted. (ii) Among all procedures where the prior probability of rejecting
$H_{0i},$ when it is true, is less than or equal to $M(RB_{i}(\psi
_{0i}\,|\,X)<q_{R}\,|\,\psi_{0i})$, then the relative belief multiple testing
algorithm maximizes the prior probability that a randomly chosen hypothesis is rejected.
\end{corollary}

\noindent\textbf{Proof}: For (i) consider a procedure for multiple testing and
let $A_{i}$ be the set of data values where $H_{0i}$ is accepted. Then by
hypothesis $M(RB_{i}(\psi_{0i}\,|\,X)>q_{A}\,|\,\psi_{0i})\leq M(A_{i}%
\,|\,\psi_{0i})$ and by the analog of Theorem \ref{thm1} $M(A_{i})\geq
M(RB_{i}(\psi_{0i}\,|\,X)$\newline$>q_{A}).$ Applying this to a randomly
chosen $H_{0i}$ gives the result. The proof of (ii) is basically the
same.\smallskip

\noindent Applying the same discussion as after Theorem \ref{thm1}, it is seen
that, under reasonable conditions, the relative belief multiple testing
algorithm minimizes the prior probability of accepting a randomly chosen
$H_{0i}$ when it is false and maximizes the prior probability of rejecting a
randomly chosen $H_{0i}$ when it is false.

Consider now the application of the relative belief multiple testing algorithm
in the previous example.\smallskip

\noindent\textbf{Example 2. }\textit{(Example 1 continued)}

In this context (\ref{exfalspos}) equals $M(RB_{i}(\mu_{0}\,|\,x)<q_{R}%
\,|\,\mu_{0})=2(1-\Phi(a_{n}(q_{R}))$ for all $i$ and so this is the value of
(\ref{exfalspos}). Therefore, $q_{R}$ is chosen to make this number suitably
small. Table \ref{biastab1} records values for (\ref{exfalspos}). From this it
is seen that for small $n$ there can be some bias against $H_{0i}$ ($q_{R}=1$)
and so the prior probability of obtaining false positives is perhaps too
large. Table \ref{biastab1} demonstrates that choosing a smaller value of
$q_{R}$ can adequately control the prior probability of false positives.%

\begin{table}[tbp] \centering
\begin{tabular}
[c]{|c|c|c|c|c|c|c|c|}\hline
$n$ & $\lambda_{0}$ & $q_{R}$ & (\ref{exfalspos}) & $n$ & $\lambda_{0}$ &
$q_{R}$ & (\ref{exfalspos})\\\hline
$1$ & $1$ & $1$ & $0.239$ $(0.228)$ & $5$ & $1$ & $1$ & $0.143$ $(0.097)$%
\\\hline
&  & $1/2$ & $0.041$ $(0.030)$ &  &  & $1/2$ & $0.051$ $(0.022)$\\\hline
&  & $1/10$ & $0.001$ $(0.000$ &  &  & $1/10$ & $0.006$ $(0.001)$\\\hline
& $2$ & $1$ & $0.156$ $(0.146)$ &  & $2$ & $1$ & $0.074$ $(0.041)$\\\hline
&  & $1/2$ & $0.053$ $(0.045)$ &  &  & $1/2$ & $0.031$ $(0.013)$\\\hline
&  & $1/10$ & $0.005$ $(0.004)$ &  &  & $1/10$ & $0.005$ $(0.001)$\\\hline
& $10$ & $1$ & $0.031$ $(0.026)$ &  & $10$ & $1$ & $0.013$ $(0.004)$\\\hline
&  & $1/2$ & $0.014$ $(0.011)$ &  &  & $1/2$ & $0.006$ $(0.002)$\\\hline
&  & $1/10$ & $0.002$ $(0.002)$ &  &  & $1/10$ & $0.001$ $(0.001)$\\\hline
\end{tabular}
\caption{Prior probability a randomly chosen hypothesis produces a  false positive,  continuous and discretized ( ) versions, in Example 2.}\label{biastab1}%
\end{table}%

For false negatives consider (\ref{exfalsnegs}) where%
\[
M(RB_{i}(\mu_{0}\,|\,x)>q_{A}\,|\,\mu_{0}\pm\delta/2)=\left\{
\begin{array}
[c]{cc}%
\begin{array}
[c]{c}%
\Phi(\sqrt{n}\delta/2\sigma+a_{n}(q_{A}))-\\
\Phi(\sqrt{n}\delta/2\sigma-a_{n}(q_{A})),
\end{array}
&
\begin{array}
[c]{c}%
1\leq q_{A}^{2}\leq\\
n\lambda_{0}^{2}+1
\end{array}
\\
0, & q_{A}^{2}>n\lambda_{0}^{2}+1.
\end{array}
\right.
\]
for all $i.$ It is easy to show that this is monotone decreasing in $\delta$
and so it is an upper bound on the expected proportion of false negatives
among those hypotheses that are actually false. The cutoff $q_{A}$ can be
chosen to make this number as small as desired. When $\delta/\sigma
\rightarrow\infty,$ then (\ref{exfalsnegs}) converges to 0 and increases to
$2\Phi(a_{n}(q_{A}))-1$ as $\delta/\sigma\rightarrow0.$ Table \ref{biastab3}
records values for (\ref{exfalsnegs}) when $\delta/\sigma=1$ so that the
$\mu_{i}$ differ from $\mu_{0}$ by one half of a standard deviation. There is
clearly some improvement but still the biasing in favor of false negatives is
readily apparent. It would seem that taking $q_{A}=\sqrt{n\lambda_{0}^{2}+1}$
gives the best results but this could be considered as quite conservative. It
is also worth remarking that all the entries in Table \ref{biastab3} can be
considered as very conservative when large effect sizes are expected.%

\begin{table}[tbp] \centering
\begin{tabular}
[c]{|c|c|c|c|c|c|c|c|}\hline
$n$ & $\lambda_{0}$ & $q_{A}$ & (\ref{exfalsnegs}) & $n$ & $\lambda_{0}$ &
$q_{A}$ & (\ref{exfalsnegs})\\\hline
$1$ & $1$ & $1.0$ & $0.704$ $(0.715)$ & $5$ & $1$ & $1.0$ & $0.631$
$(0.702)$\\\hline
&  & $1.2$ & $0.527$ $(0.503)$ &  &  & $2.0$ & $0.302$ $(0.112)$\\\hline
&  & $1.4$ & $0.141$ $(0.000)$ &  &  & $2.4$ & $0.095$ $(0.000)$\\\hline
& $2$ & $1.0$ & $0.793$ $(0.805)$ &  & $2$ & $1.0$ & $0.747$ $(0.822)$\\\hline
&  & $2.0$ & $0.359$ $(0.304)$ &  &  & $3.0$ & $0.411$ $(0.380)$\\\hline
&  & $2.2$ & $0.141$ $(0.000)$ &  &  & $4.5$ & $0.084$ $(0.000)$\\\hline
& $10$ & $1.0$ & $0.948$ $(0.955)$ &  & $10$ & $1.0$ & $0.916$ $(0.961)$%
\\\hline
&  & $5.0$ & $0.708$ $(0.713)$ &  &  & $10.0$ & $0.552$ $(0.588)$\\\hline
&  & $10.0$ & $0.070$\ $(0.000)$ &  &  & $22.0$ & $0.080$ $(0.000)$\\\hline
\end{tabular}
\caption{Prior probability a randomly chosen hypothesis produces a  false negative when $\delta/\sigma =1$,  continuous and discretized ( ) versions, in Example 2.}\label{biastab3}%
\end{table}%

Now consider the situation that led to Figure 2. For this $k=1000,n=5$ and
$\lambda_{0}=1.94$ is the elicited value. From Table \ref{biastab1} with
$q_{R}=1.0,$ about 8\% false positives are expected \textit{a priori} and from
Table \ref{biastab3} with $q_{A}=1.0,$ a worst case upper bound on the
\textit{a priori} expected percentage of false negatives is about 75\%. The
top part of Table \ref{confusion1} indicates that with $q_{R}=q_{A}=1.0,$ then
4.9\% (34 of 700) false positives and 0.1\% (3 of 300) false negatives were
obtained. With these choices of the cutoffs all hypotheses are classified.
Certainly the upper bound 75\% seems far too pessimistic in light of the
results, but recall that Table \ref{biastab3} is computed at the false values
$\mu=\pm0.5.$ The relevant \textit{a priori} expected percentage of false
negatives when $\mu=\pm2.0$ is about 3.5\%. The bottom part of Table
\ref{confusion1} gives the relevant values when $q_{R}=0.5$ and $q_{A}=3.0.$
In this case there are 2.1\% (9 of 428) false positives and 0\% false
negatives but 39.9\% (272 out of 700) of the true hypotheses and 4.3\% (13 out
300)\ of the false hypotheses were not classified as the relevant relative
belief ratio lay between $q_{R}$ and $q_{A}.$ So in this case being more
conservative has reduced the error rates with the price being a large
proportion of the true hypotheses don't get classified. The procedure has
worked well in this example but of course the error rates can be expected to
rise when the false values move towards the null and improve when they move
away from the null. $\blacksquare$%

\begin{table}[tbp] \centering
\begin{tabular}
[c]{|l|l|l|}\hline
Decision & $\mu=0$ & $\mu=2$\\\hline
Accept $\mu=0$ using $q_{A}=1.0$ & \multicolumn{1}{|c|}{666} &
\multicolumn{1}{|c|}{3}\\
Reject $\mu=0$ using $q_{R}=1.0$ & \multicolumn{1}{|c|}{34} &
\multicolumn{1}{|c|}{297}\\
Not classified & \multicolumn{1}{|c|}{0} & \multicolumn{1}{|c|}{0}\\\hline
Accept $\mu=0$ using $q_{A}=3.0$ & \multicolumn{1}{|c|}{419} &
\multicolumn{1}{|c|}{0}\\
Reject $\mu=0$ using $q_{R}=0.5$ & \multicolumn{1}{|c|}{9} &
\multicolumn{1}{|c|}{287}\\
Not classified & \multicolumn{1}{|c|}{272} & \multicolumn{1}{|c|}{13}\\\hline
\end{tabular}
\caption{Confusion matrices for Example 2 with $k=1000$ when 700 of the $\mu_i$ equal 0 and 300 equal 2.
}\label{confusion1}%
\end{table}%

What is implemented in an application depends on the goals. If the primary
purpose is to protect against false positives, then Table \ref{biastab1}
indicates that this is accomplished fairly easily. Protecting against false
negatives is more difficult. Since the actual effect sizes are not known a
decision has to be made. Note that choosing a cutoff is equivalent to saying
that one will only accept $H_{0i}$ if the belief in the truth of $H_{0i}$ has
increased by a factor at least as large as $q_{A}.$ Computations such as in
Table \ref{biastab3} can be used to provide guidance but there is no avoiding
the need to be clear about what effect sizes are deemed to be important or the
need to obtain more data when this is necessary. One comforting aspect of the
methodology is that error rates are effectively controlled but there may be
many true hypotheses not classified.

The idea of controlling the prior probability of a randomly chosen hypothesis
yielding a false positive or a false negative via (\ref{exfalspos}) or
(\ref{exfalsnegs}), respectively, can be extended. For example, consider the
prior probability that a random sample of $l$ from $k$ hypotheses yields at
least one false positive%
\begin{equation}
\frac{1}{\binom{k}{l}}\sum_{\{i_{1},\ldots,i_{l}\}\subset\{1,\ldots
,k\}}M\left(
\begin{array}
[c]{c}%
\text{at least one of }RB_{i_{j}}(\psi_{0i_{j}}\,|\,X)<q_{R}\\
\text{for }j=1,\ldots,l\,|\,\psi_{0i_{1}},\ldots,\psi_{0i_{l}}%
\end{array}
\right)  . \label{lrandom}%
\end{equation}
In the context of the examples of this paper, and many others, the term in
(\ref{lrandom}) corresponding to $\{i_{1},\ldots,i_{l}\}$ equals $M($at least
one of $RB_{i_{j}}(\psi_{0i_{j}}\,|\,X)<q_{R}$ for $j=1,\ldots,l\,|\,\psi
_{0}).$ The following result, whose proof is given in the Appendix, then leads
to an interesting property for (\ref{lrandom}).

\begin{lemma}
\label{combin}For probability model $(\Omega,\mathcal{F},P),$ the probability
that at least one of $l\leq k$ randomly selected events from $\{A_{1}%
,\ldots,A_{k}\}\subset\mathcal{F}$ occurs is an upper bound on the probability
that at least one of $l^{\prime}\leq l$ randomly selected events from
$\{A_{1},\ldots,A_{k}\}\subset\mathcal{F}$ occur.
\end{lemma}

\noindent It then follows, by taking $A_{i}=\{x:RB_{i}(\psi_{0i}%
\,|\,x)<q_{R}\},$ that (\ref{lrandom}) is an upper bound on the prior
probability that a random sample of $l^{\prime}$ hypotheses yields at least
one false positive whenever $l^{\prime}\leq l.$ So (\ref{lrandom}) leads to a
more rigorous control over the possibility false positives, if so desired. A
similar result is obtained for false negatives.

\section{Applications}

One application of the relative belief multiple testing algorithm is to the
problem of inducing sparsity.\smallskip

\noindent\textbf{Example 3.} \textit{Testing for sparsity.}

The position taken here is that, rather than trying to induce sparsity through
an estimation procedure, an assessment of sparsity is made through an explicit
measure of the evidence as to whether or not a particular assignment is valid.
Virtues of this approach are that it is based on a measure of evidence and it
is not dependent on the form of the prior as any prior can be used. Also, it
has some additional advantages when prior-data conflict is taken into account.

Consider the context of Example 1. A natural approach to inducing sparsity is
to estimate $\mu_{i}$ by $\mu_{0}$ whenever $RB_{i}(\mu_{0}\,|\,x)>q_{A}.$
From the simulation it is seen that this works extremely well when $q_{A}=1$
for both $k=10$ and $k=1000.$ It also works when $k=1000$ and $q_{A}=3$, in
the sense that the error rate is low, but it is conservative in the amount of
sparsity it induces in that case. Again the goals of the application will
dictate what is appropriate.

A common method for inducing sparsity is to use a penalized estimator as in
the LASSO introduced in Tibshirani (1996). It is noted there that the
LASSO\ is equivalent to MAP\ estimation when using a product of independent
Laplace priors. This aspect was pursued further in Park and Casella (2006)
which adopted a more formal Bayesian approach but still used MAP.

Consider then a product of independent Laplace priors for the prior on $\mu,$
namely, $(\sqrt{2}\lambda_{0}\sigma)^{-k}\exp\{-(\sqrt{2}/\lambda_{0}%
\sigma)\sum\nolimits_{i=1}^{k}|\mu_{i}-\mu_{0}|\}$ where $\sigma$ is assumed
known and $\mu_{0},\lambda_{0}$ are hyperparameters. Note that each Laplace
prior has mean $\mu_{0}$ and variance $\lambda_{0}^{2}\sigma^{2}.$ Using the
elicitation algorithm provided in Example 1, but replacing the normal prior
with a Laplace prior, leads to the assignment $\mu_{0}=$ $(m_{l}+m_{u})/2,$
$\lambda_{0}=(m_{u}-m_{l})/2\sigma G^{-1}(0.995)$ where $G^{-1}(p)=2^{-1/2}%
\log2p$ when $p\leq1/2,G^{-1}(p)=-2^{-1/2}\log2(1-p)$ when $p\geq1/2\ $and
$G^{-1}$ denotes the quantile function of a Laplace distribution with mean $0$
and variance $1.$ With the specifications used in the simulations of Example
1, this leads to $\mu_{0}=0$ and $\lambda_{0}=1.54$ which implies a smaller
variance than the value $\lambda_{0}=1.94$ used with the normal prior and so
the Laplace prior is more concentrated about 0.

The posteriors for the $\mu_{i}$ are independent with the density for $\mu
_{i}$ proportional to $\exp\{-n(\bar{x}_{i}-\mu_{i})^{2}/2\sigma^{2}-\sqrt
{2}|\mu_{i}-\mu_{0}|/\lambda_{0}\sigma\}$ giving the MAP estimator%
\[
\mu_{i\text{MAP}}(x)=\left\{
\begin{array}
[c]{llc}%
\bar{x}_{i}+\sqrt{2}\sigma/\lambda_{0}n, &  & \bar{x}_{i}<\mu_{0}-\sqrt
{2}\sigma/\lambda_{0}n\\
\mu_{0}, &  & \mu_{0}-\sqrt{2}\sigma/\lambda_{0}n\leq\bar{x}_{i}\leq\mu
_{0}+\sqrt{2}\sigma/\lambda_{0}n\\
\bar{x}_{i}-\sqrt{2}\sigma/\lambda_{0}n, &  & \bar{x}_{i}>\mu_{0}+\sqrt
{2}\sigma/\lambda_{0}n.
\end{array}
\right.
\]
The MAP estimate of $\mu_{i}$ is sometimes forced to equal $\mu_{0}$ although
this effect is negligible whenever $\sqrt{2}\sigma/\lambda_{0}n$ is small.

The LASSO induces sparsity through estimation by taking $\lambda_{0}$ to be
small. By contrast the evidential approach, based on the normal prior and the
relative belief ratio, induces sparsity through taking $\lambda_{0}$ large.
The advantage to this latter approach is that by taking $\lambda_{0}$ large,
prior-data conflict is avoided. When taking $\lambda_{0}$ small, the potential
for prior-data conflict rises as the true values can be deep into the tails of
the prior. For example, for the simulations of Example 1 $\sqrt{2}%
\sigma/\lambda_{0}n=0.183$ which is smaller than the $\delta/2=0.5$ used in
the relative belief approach with the normal prior.\ So it can be expected
that the LASSO will do worse here and this is reflected in Table
\ref{confusion2} where there are far too many false negatives. To improve this
the value of $\lambda_{0}$ needs to be reduced although note that this is
determined by an elicitation and there is the risk of then encountering
prior-data conflict. Another possibility is to implement the evidential
approach with the elicited Laplace prior and the discretization as done with
the normal prior and then we can expect results similar to those obtained in
Example 1.%

\begin{table}[tbp] \centering
\begin{tabular}
[c]{|l|l|l|}\hline
Decision & $\mu=0$ & $\mu=2$\\\hline
Accept $\mu=0$ using $q_{A}=1.0$ & \multicolumn{1}{|c|}{227} &
\multicolumn{1}{|c|}{0}\\
Reject $\mu=0$ using $q_{A}=1.0$ & \multicolumn{1}{|c|}{473} &
\multicolumn{1}{|c|}{300}\\\hline
\end{tabular}
\caption{Confusion matrices using LASSO with $k=1000$ when 700 of the $\mu_i$ equal 0 and 300 equal 2 in Example 3
}\label{confusion2}
\end{table}%

It is also interesting to compare the MAP\ estimation approach and the
relative belief approach with respect to the conditional prior probabilities
of $\mu_{i}$ being assigned the value $\mu_{0}$ when the true value actually
is $\mu_{0}.$ It is easily seen that, based on the Laplace prior,
$M(\mu_{i\text{MAP}}(x)=\mu_{0}\,|\,\mu_{0})=2\Phi(\sqrt{2}/\lambda_{0}%
\sqrt{n})-1$ and this converges to 0 as $n\rightarrow\infty$ or $\lambda
_{0}\rightarrow\infty.$ For the relative belief approach $M(RB_{i}(\mu
_{0}\,|\,x)>q_{A}\,|\,\mu_{0})$ is the relevant probability. With either the
normal or Laplace prior $M(RB_{i}(\mu_{0}\,|\,x)>q_{A}\,|\,\mu_{0})$ converges
to 1 both as $n\rightarrow\infty$ and as $\lambda_{0}\rightarrow\infty.$ In
particular, with enough data the correct assignment is always made using
relative belief but not with MAP based on the Laplace prior.

While the Laplace and normal priors work equally with the relative belief
multiple testing algorithm, there don't appear to be any advantages to using
the Laplace prior. One could argue too that the singularity of the Laplace
prior at its mode makes it an odd choice and there doesn't seem to be a good
justification for this.\ Furthermore, the computations are harder with the
Laplace prior, particularly with more complex models. So using a normal prior
seems preferable overall. $\blacksquare\smallskip$

An example with considerable practical significance is now
considered.$\smallskip$

\noindent\textbf{Example 4.} \textit{Full rank regression.}

Suppose the basic model is given by $y=\beta_{0}+\beta_{1}x_{1}+\cdots
+\beta_{k}x_{k}+z=\beta_{0}+x^{\prime}\beta_{1:k}+z$ where the $x_{i}$ are
predictor variables, $z\sim N(0,\sigma^{2})$ and the $\beta_{i}$ and
$\sigma^{2}$ are unknown. The main issue in this problem is testing
$H_{0i}:\beta_{i}=0$ for $i=1,\ldots,k$ to establish which variables have any
effect on the response. The prior distribution of $(\beta,\sigma^{2})$ is
taken to be
\begin{equation}
\beta\,|\,\sigma^{2}\sim N_{k+1}(0,\sigma^{2}\Sigma_{0}),1/\sigma^{2}%
\sim\text{gamma}_{\text{rate}}(\alpha_{1},\alpha_{2}), \label{regprior}%
\end{equation}
for some hyperparameters $\Sigma_{0}$ and $(\alpha_{1},\alpha_{2}).$ Note that
this may entail subtracting a known, fixed constant from each\ $y$ value so
that the prior for $\beta_{0}$ is centered at 0. Taking 0 as the central value
for the priors on the remaining $\beta_{i}$ seems appropriate when the primary
concern is whether or not each $x_{i}$ is having any effect. Also, it will be
assumed that the observed values of the predictor variable have been
standardized so that for observations $(y,X)\in R^{n}\times R^{n\times(k+1)},$
where $X=(\mathbf{1},\mathbf{x}_{1},\ldots,\mathbf{x}_{k}),$ then
$\mathbf{1}^{\prime}\mathbf{x}_{i}=0$ and $||\mathbf{x}_{i}||^{2}=1$ for
$i=1,\ldots,k.$ The marginal prior for $\beta_{i}$ is then $\{(\alpha
_{2}/\alpha_{1})\sigma_{0ii}^{2}\}^{1/2}t_{2\alpha_{1}}$ where $t_{2\alpha
_{1}}$ denotes the $t$ distribution on $2\alpha_{1}$ degrees of freedom, for
$i=0,\ldots,k.$ Hereafter, we will take $\Sigma_{0}=\lambda_{0}^{2}I_{k+1}$
although it is easy to generalize to more complicated choices.

The elicitation of the hyperparameters is carried out via an extension of a
method developed in Cao, Evans and Guttman (2014) for the multivariate normal
distribution. Suppose that it is known with virtual certainty, based on our
knowledge of the measurements being taken, that $\beta_{0}+x^{\prime}%
\beta_{1:k}$ will lie in the interval $(-m_{0},m_{0})$ for some $m_{0}>0$ for
all $x\in R$ where $R$ is a compact set centered at 0. On account of the
standardization, $R\subset\lbrack-1,1]^{k}.$ Again `virtual certainty' is
interpreted as probability greater than or equal to $\gamma$ where $\gamma$ is
some large probability like $0.99.$ Therefore, the prior on $\beta$ must
satisfy $2\Phi(m_{0}/\sigma\lambda_{0}\{1+x^{\prime}x\}^{1/2})-1\geq\gamma$
for all $x\in R$ and this implies that%
\begin{equation}
\sigma\leq m_{0}/\lambda_{0}\tau_{0}z_{(1+\gamma)/2} \label{ineq1}%
\end{equation}
where $\tau_{0}^{2}=1+\max_{x\in R}||x||^{2}\leq1+k$ with equality when
$R=[-1,1]^{k}.$

An interval that will contain a response value $y$ with virtual certainty,
given predictor values $x,$ is $\beta_{0}+x^{\prime}\beta_{1:k}\pm\sigma
z_{(1+\gamma)/2}.$ Suppose that we have lower and upper bounds $s_{1}$ and
$s_{2}$ on the half-length of this interval so that $s_{1}\leq\sigma
z_{(1+\gamma)/2}\leq s_{2}$ or, equivalently,%
\begin{equation}
s_{1}/z_{(1+\gamma)/2}\leq\sigma\leq s_{2}/z_{(1+\gamma)/2} \label{ineq2}%
\end{equation}
holds with virtual certainty. Combining (\ref{ineq2}) with (\ref{ineq1})
implies $\lambda_{0}=m_{0}/s_{2}\tau_{0}.$

To obtain the relevant values of $\alpha_{1}$ and $\alpha_{2}$ let $G\left(
\alpha_{1},\alpha_{2},\cdot\right)  $ denote the cdf of the
gamma$_{\text{rate}}\left(  \alpha_{1},\alpha_{2}\right)  $ distribution and
note that $G\left(  \alpha_{1},\alpha_{2},z\right)  =G\left(  \alpha
_{1},1,\alpha_{2}z\right)  .$ Therefore, the interval for $1/\sigma^{2}$
implied by (\ref{ineq2}) contains $1/\sigma^{2}$ with virtual certainty, when
$\alpha_{1},\alpha_{2}$ satisfy $G^{-1}(\alpha_{1},\alpha_{2},(1+\gamma
)/2)=s_{1}^{-2}z_{(1+\gamma)/2}^{2},G^{-1}(\alpha_{1},\alpha_{2}%
,(1-\gamma)/2)=s_{2}^{-2}z_{(1-\gamma)/2}^{2},$ or equivalently
\begin{align}
G(\alpha_{1},1,\alpha_{2}s_{1}^{-2}z_{(1+\gamma)/2}^{2})  &  =(1+\gamma
)/2,\label{eq3}\\
G(\alpha_{1},1,\alpha_{2}s_{2}^{-2}z_{(1-\gamma)/2}^{2})  &  =(1-\gamma)/2.
\label{eq4}%
\end{align}
It is a simple matter to solve these equations for $\left(  \alpha_{1}%
,\alpha_{2}\right)  .$ For this choose an initial value for $\alpha_{1}$ and,
using (\ref{eq3}), find $z$ such that $G(\alpha_{1},1,z)=(1+\gamma)/2,$ which
implies $\alpha_{2}=z/s_{1}^{-2}z_{(1+\gamma)/2}^{2}.$ If the left-side of
(\ref{eq4}) is less (greater) than $(1-\gamma)/2,$ then decrease
(increase)\ the value of $\alpha_{1}$ and repeat step 1. Continue iterating
this process until satisfactory convergence is attained.

The methods discussed in Evans and Moshonov (2006) are available for checking
the prior to see if it is contradicted by the data. Methods are specified
there for checking each of the components in the hierarchy, namely, first
check the prior on $\sigma^{2}$ and, if it passes, then check the prior on
$\beta.$ If conflict is found, then the methods discussed in Evans and Jang
(2011b) are available to modify the prior appropriately.

Assuming that $X$ is of rank $k+1,$ the posterior of $(\beta,\sigma^{2})$ is
given by
\begin{align}
\beta\,|\,y,\sigma^{2}  &  \sim N_{k+1}(\beta(X,y),\sigma^{2}\Sigma
(X)),\nonumber\\
1/\sigma^{2}\,|\,y  &  \sim\text{gamma}_{\text{rate}}((n+2\alpha_{1}%
)/2,\alpha_{2}(X,y)/2), \label{regpost}%
\end{align}
where $b=(X^{\prime}X)^{-1}X^{\prime}y,\beta(X,y)=\Sigma(X)X^{\prime}%
Xb,\Sigma(X)=(X^{\prime}X+\Sigma_{0}^{-1})^{-1}$ and $\alpha_{2}%
(X,y)=||y-Xb||^{2}+(Xb)^{\prime}(I_{n}-X\Sigma(X)X^{\prime})Xb+2\alpha_{2}.$
Then the marginal posterior for $\beta_{i}$ is given by $\beta_{i}%
(X,y)+\{\alpha_{2}(X,y)\sigma_{ii}(X)/(n+2\alpha_{1})\}^{1/2}t_{n+2\alpha_{1}%
}$ and the relative belief ratio for $\beta_{i}$ at $0$ equals%
\begin{align}
RB_{_{i}}(0\,|\,X,y)=  &  \frac{\Gamma\left(  \frac{n+2\alpha_{1}+1}%
{2}\right)  \Gamma\left(  \alpha_{1}\right)  }{\Gamma\left(  \frac{2\alpha
_{1}+1}{2}\right)  \Gamma\left(  \frac{n+2\alpha_{1}}{2}\right)  }\left(
1+\frac{\beta_{i}^{2}(X,y)}{\alpha_{2}(X,y)\sigma_{ii}(X)}\right)
^{-\frac{n+2\alpha_{1}+1}{2}}\times\nonumber\\
&  \left(  \frac{\alpha_{2}(X,y)\sigma_{ii}(X)}{\alpha_{2}^{2}\lambda_{0}^{2}%
}\right)  ^{-\frac{1}{2}}. \label{rbreg}%
\end{align}

Rather than using (\ref{rbreg}), however, the distributional results are used
to compute the discretized relative belief ratios as in Example 1. For this
$\delta>0$ is required to determine an appropriate discretization and it will
be assumed here that this is the same for all the $\beta_{i},$ although the
procedure can be easily modified if this is not the case in practice. Note
that such a $\delta$ is effectively determined by the amount that $x_{i}%
\beta_{i}$ will vary from $0$ for $x\in R.$ Since $x_{i}\in\lbrack-1,1]$ then
$|x_{i}\beta_{i}|\leq\delta$ provided $|\beta_{i}|\leq\delta.$ When this
variation is suitably small as to be immaterial, then such a $\delta$ is
appropriate for saying $\beta_{i}$ is effectively $0.$ Note that determination
of the hyperparameters and $\delta$ is dependent on the application.

Again inference can be made concerning $\xi=\Xi(\beta,\sigma^{2}),$ the
proportion of the $\beta_{i}$ effectively equal to $0.$ As in Example 1,
however, we can expect bias when the amount of variability in the data is
large relative to $\delta$ or the prior is too diffuse. To implement the
relative belief multiple testing algorithm the quantities (\ref{exfalspos})
and (\ref{exfalsnegs}) need to be computed to determine $q_{R}$ and $q_{A},$
respectively. The conditional prior distribution of $(b,||y-Xb||^{2})$, given
$(\beta,\sigma^{2}),$ is $b\sim N_{k+1}(\beta,\sigma^{2}(X^{\prime}X)^{-1})$
statistically independent of $||y-Xb||^{2}\sim$ gamma$((n-k-1)/2,\sigma
^{-2}/2).$ So computing (\ref{exfalspos}) and (\ref{exfalsnegs}) can be
carried out by generating $(\beta,\sigma^{2})$ from the relevant conditional
prior, generating $(b,||y-Xb||^{2})$ given $(\beta,\sigma^{2}),$ and using
(\ref{rbreg}).

To illustrate these computations the diabetes data set discussed in Efron,
Hastie, Johnstone and Tibshirani (2006) and Park and Casella (2008) is now
analyzed. With $\gamma=0.99,$ the values $m_{0}=100,s_{1}=75,s_{2}=200$ were
used to determine the prior together with $\tau_{0}=1.05$ determined from the
$X$ matrix. This lead to the values $\lambda_{0}=0.48,\alpha_{1}%
=7.29,\alpha_{2}=13641.35$ being chosen for the hyperparameters. Using the
methods developed in$\ $Evans and Moshonov (2006), a first check was made on
the prior on $\sigma^{2}$ against the data and a tail probability equal to
$0.19$ was obtained indicating there is no prior-data conflict with this
prior. Given no prior-data conflict at the first stage, the prior on $\beta$
was then checked and the relevant tail probability of 0.00 was obtained
indicating a strong degree of conflict. Following the argument in Evans and
Jang (2011) the value of $\lambda_{0}$ was increased to choose a prior weakly
informative with respect to our initial choice and this lead to choosing the
value $\lambda_{0}=5.00$ and then the relevant tail probability equals $0.32.$

Using this prior, the relative belief estimates, ratios and strengths are
recorded in Table \ref{fullrankregests}. From this it is seen that there is
strong evidence against $\beta_{i}=0$ for the variables sex, bmi, map and ltg
and no evidence against $\beta_{i}=0$ for any other variables. There is strong
evidence of in favor of $\beta_{i}=0$ for age and ldl, moderate evidence in
favor of $\beta_{i}=0$ for the constant, tc, tch and glu and perhaps only weak
evidence in favor of $\beta_{i}=0$ for hdl.

As previously discussed it is necessary to consider the issue of bias, namely,
compute the prior probability of getting a false positive for different
choices of $q_{R}$ and the prior probability of getting a false negative for
different choices of $q_{A}.$ The value of (\ref{exfalspos}) is $0.0003$ when
$q_{R}=1$ and so there is virtually no bias in favor of false positives and
one can feel confident that the predictors identified as having an effect do
so. The story is somewhat different, however, when considering the possibility
of false negatives via (\ref{exfalsnegs}). For example, with $q_{A}=1,$ then
(\ref{exfalsnegs}) equals $0.9996$ and when $q_{A}=100$ then (\ref{exfalsnegs}%
) equals $0.7998$. So there is substantial bias in favor of the null
hypotheses and undoubtedly this is due to the diffuseness of the prior. The
implication is that we cannot be entirely confident concerning those
$\beta_{i}$ assigned to be equal to 0. Recall, that the first prior proposed
lead to prior-data conflict and as such a much more diffuse prior was
substituted. The bias in favor of false negatives could be mitigated by making
the prior less diffuse. It is to be noted, however, that this is an exercise
that should be conducted prior to collecting the data as there is a danger
that the choice of the prior will be too heavily influenced by the observed
data. The real cure for any bias in an application is to collect more data.
$\blacksquare\smallskip$%

\begin{table}[tbp] \centering
\begin{tabular}
[c]{|c|c|c|c|}\hline
Variable & Estimates & $RB_{_{i}}(0\,|\,X,y)$ & Strength\\\hline
Constant & \multicolumn{1}{|r|}{$2$} & \multicolumn{1}{|r|}{$2454.86$} &
$0.44$\\\hline
age & \multicolumn{1}{|r|}{$-4$} & \multicolumn{1}{|r|}{$153.62$} &
$0.95$\\\hline
sex & \multicolumn{1}{|r|}{$-224$} & \multicolumn{1}{|r|}{$0.13$} &
$0.00$\\\hline
bmi & \multicolumn{1}{|r|}{$511$} & \multicolumn{1}{|r|}{$0.00$} &
$0.00$\\\hline
map & \multicolumn{1}{|r|}{$314$} & \multicolumn{1}{|r|}{$0.00$} &
$0.00$\\\hline
tc & \multicolumn{1}{|r|}{$162$} & \multicolumn{1}{|r|}{$33.23$} &
$0.36$\\\hline
ldl & \multicolumn{1}{|r|}{$-20$} & \multicolumn{1}{|r|}{$57.65$} &
$0.90$\\\hline
hdl & \multicolumn{1}{|r|}{$167$} & \multicolumn{1}{|r|}{$27.53$} &
$0.15$\\\hline
tch & \multicolumn{1}{|r|}{$114$} & \multicolumn{1}{|r|}{$49.97$} &
$0.37$\\\hline
ltg & \multicolumn{1}{|r|}{$496$} & \multicolumn{1}{|r|}{$0.00$} &
$0.00$\\\hline
glu & \multicolumn{1}{|r|}{$77$} & \multicolumn{1}{|r|}{$66.81$} &
$0.23$\\\hline
\end{tabular}
\caption{Relative belief estimates, relative belief ratios and strengths for assessing no effect  for the diabetes data in Example 4..}\label{fullrankregests}%
\end{table}%

Next we consider the application to regression with $k+1>n.$\smallskip

\noindent\textbf{Example 5.} \textit{Non-full rank regression.}

In a number of applications $k+1>n$ and so $X$ is of rank $l<n$. In this
situation, suppose $\{\mathbf{x}_{1},\ldots,\mathbf{x}_{l}\}$ forms a basis
for $\mathcal{L}(\mathbf{x}_{1},\ldots,\mathbf{x}_{k}),$ perhaps after
relabeling the predictors, and write $X=(\mathbf{1}$ $X_{1}$ $X_{2})$\textbf{
}where $X_{1}=(\mathbf{x}_{1}\,\ldots\,\mathbf{x}_{l}).$ For given
$r=(X_{1}\,X_{2})\beta_{1:k}$ there will be many solutions $\beta_{1:k}.$ A
particular solution is given by $\beta_{1:k\ast}=(X_{1}(X_{1}^{\prime}%
X_{1})^{-1}$ $0)^{\prime}r.$ The set of all solutions is then given by
$\beta_{1:k\ast}+\ker(X_{1}\,X_{2})$ where $\ker(X_{1}\,X_{2})=\{(-B^{\prime}$
$I_{k-l})^{\prime}\eta:\eta\in R^{k-l}\},B=(X_{1}^{\prime}X_{1})^{-1}%
X_{1}^{\prime}X_{2}$ and the columns of $C=(-B^{\prime}$ $I_{k-l})^{\prime}$
give a basis for $\ker(X_{1}\,X_{2}).$ Given that sparsity is expected for the
true $\beta_{1:k},$ it is natural to consider the solution which minimizes
$||\beta_{1:k}||^{2}$ for $\beta_{1:k}\in\beta_{1:k\ast}+\mathcal{L}(C).$
Using $\beta_{1:k\ast}$, and applying the Sherman-Morrison-Woodbury formula to
$C(C^{\prime}C)^{-1}C^{\prime},$ this is given by the Moore-Penrose solution
\begin{equation}
\beta_{1:k}^{MP}=(I_{k}-C(C^{\prime}C)^{-1}C^{\prime})\beta_{1:k\ast}%
=(I_{l}\text{ }B)^{\prime}\omega_{1:l} \label{sparsesolution}%
\end{equation}
where $\omega_{1:l}=(I_{l}+BB^{\prime})^{-1}(\beta_{1:l}+B\beta_{l+1:k}).$

From (\ref{regprior}) with $\Sigma_{0}=\lambda_{0}^{2}I_{k+1}$, the
conditional prior distribution of $(\beta_{0},\omega_{1:l})$ given $\sigma
^{2}$ is $\beta_{0}\,|\,\sigma^{2}\sim N(0,\sigma^{2}\lambda_{0}^{2})$
independent of $\omega_{1:l}\,|\,\sigma^{2}\sim N_{l}(0,\sigma^{2}\lambda
_{0}^{2}(I_{l}+BB^{\prime})^{-1})$ which, using (\ref{sparsesolution}),
implies $\beta_{1:k}^{MP}\,|\,\sigma^{2}\sim N_{k}(0,\sigma^{2}\Sigma
_{0}(B)),$ conditionally independent of $\beta_{0},$ where
\[
\Sigma_{0}(B)=\lambda_{0}^{2}\left(
\begin{array}
[c]{cc}%
(I_{l}+BB^{\prime})^{-1} & (I_{l}+BB^{\prime})^{-1}B\\
B^{\prime}(I_{l}+BB^{\prime})^{-1} & B^{\prime}(I_{l}+BB^{\prime})^{-1}B
\end{array}
\right)  .
\]
With $1/\sigma^{2}\sim\,$gamma$_{\text{rate}}(\alpha_{1},\alpha_{2}),$ this
implies that the unconditional prior of the $i$-th coordinate of $\beta
_{1:k}^{MP}$ is $\left(  \lambda_{0}^{2}\alpha_{2}\sigma_{ii}^{2}%
(B)/\alpha_{1}\right)  ^{1/2}t_{2\alpha_{1}}.$

Putting $X_{\ast}=(\mathbf{1}$ $X_{1}+X_{2}B^{\prime})$ gives the full rank
model $y\,|\,\beta_{0},\omega_{1:l},\sigma^{2}\sim N_{n}(X_{\ast}(\beta
_{0},\omega_{1:l}^{\prime})^{\prime},\sigma^{2}I_{n}).$ As in Example 4 then
$(\beta_{0},\omega_{1:l})\,|\,y,\sigma^{2}\sim N_{l}(\omega(X_{\ast}%
,y),$\newline$\sigma^{2}\Sigma(X_{\ast})),1/\sigma^{2}\,|\,y\sim
\,$gamma$_{\text{rate}}((n+2\alpha_{1})/2,\alpha_{2}(X_{\ast},y)/2)$ where
$\omega(X_{\ast},y)=\Sigma(X_{\ast})X_{\ast}^{\prime}X_{\ast}b_{\ast},b_{\ast
}=(X_{\ast}^{\prime}X_{\ast})^{-1}X_{\ast}^{\prime}y$ and
\begin{align*}
\Sigma^{-1}(X_{\ast})  &  =\left(
\begin{array}
[c]{cc}%
n & 0\\
0 & (X_{1}+X_{2}B^{\prime})^{\prime}(X_{1}+X_{2}B^{\prime})
\end{array}
\right)  +\lambda_{0}^{-2}\left(
\begin{array}
[c]{cc}%
1 & 0\\
0 & (I_{l}+BB^{\prime})
\end{array}
\right)  ,\\
\alpha_{2}(X_{\ast},y)  &  =||y-X_{\ast}b_{\ast}||^{2}+(X_{\ast}b_{\ast
})^{\prime}(I_{n}-X_{\ast}\Sigma(X_{\ast})X_{\ast}^{\prime})X_{\ast}b_{\ast
}+2\alpha_{2}.
\end{align*}
Now noting that $(X_{1}+X_{2}B^{\prime})^{\prime}(X_{1}+X_{2}B^{\prime
})=(I_{l}+BB^{\prime})X_{1}^{\prime}X_{1}(I_{l}+BB^{\prime}),$ this implies
$b_{\ast}^{\prime}=(\bar{y},(I_{l}+BB^{\prime})^{-1}b_{1}),$ where
$b_{1}=(X_{1}^{\prime}X_{1})^{-1}X_{1}^{\prime}y$ is the least-squares
estimate of $\beta_{1:l},$ and
\begin{align*}
\Sigma(X_{\ast})  &  =\left(
\begin{array}
[c]{cc}%
n+\lambda_{0}^{-2} & 0\\
0 & (I_{l}+BB^{\prime})X_{1}^{\prime}X_{1}(I_{l}+BB^{\prime})+\lambda_{0}%
^{-2}(I_{l}+BB^{\prime})
\end{array}
\right)  ^{-1},\\
\omega(X_{\ast},y)  &  =\Sigma(X_{\ast})X_{\ast}^{\prime}X_{\ast}b_{\ast
}=\left(
\begin{array}
[c]{c}%
n\bar{y}/(n+\lambda_{0}^{-2})\\
(I_{l}+BB^{\prime}+\lambda_{0}^{-2}(X_{1}^{\prime}X_{1})^{-1})^{-1}b_{1}%
\end{array}
\right)  .
\end{align*}
Using (\ref{sparsesolution}), then $\beta_{0}\,|\,y,\sigma^{2}\sim
N(n(n+\lambda_{0}^{-2})^{-1}\bar{y},\sigma^{2}(n+\lambda_{0}^{-2})^{-1})$
independent of $\beta_{1:k}^{MP}\,|\,y,\sigma^{2}\sim N_{k}(\beta
^{MP}(X,y),\sigma^{2}\Sigma^{MP}(X))$ where
\[
\beta^{MP}(X,y)=\left(
\begin{array}
[c]{c}%
Db_{1}\\
B^{\prime}Db_{1}%
\end{array}
\right)  ,\text{ \ \ }\Sigma^{MP}(X)=\left(
\begin{array}
[c]{cc}%
E & EB\\
B^{\prime}E & B^{\prime}EB
\end{array}
\right)
\]
with $D=(I_{l}+BB^{\prime}+\lambda_{0}^{-2}(X_{1}^{\prime}X_{1})^{-1})^{-1}$
and $E=((I_{l}+BB^{\prime})(X_{1}^{\prime}X_{1})(I_{l}+BB^{\prime}%
)+\lambda_{0}^{-2}(I_{l}+BB^{\prime}))^{-1}.$ The marginal posterior for
$\beta_{i}^{MP}$ is then given by $\beta_{i}^{MP}(X,y)+\{\alpha_{2}(X_{\ast
},y)\sigma_{ii}^{MP}(X)/(n+2\alpha_{1})\}^{1/2}t_{n+2\alpha_{1}}.$ Relative
belief inferences for the coordinates of $\beta_{1:k}^{MP}$ can now be
implemented just as in Example 4.

We consider a numerical example where there is considerable sparsity. For this
let $X_{1}\in R^{n\times l}$ be formed by taking the second through $l$-th
columns of the $(l+1)$-dimensional Helmert matrix, repeating each row $m$
times and then normalizing. So $n=m(l+1)$ and the columns of $X_{1}$ are
orthonormal and orthogonal to $\mathbf{1.}$ It is supposed that the first
$l_{1}\leq l$ of the variables giving rise to the columns of $X_{1}$ have
$\beta_{i}\neq0$ whereas the last $l-l_{1}$ have $\beta_{i}=0$ and that the
variables corresponding to the first $l_{2}\leq k-l$ columns of $X_{2}%
=X_{1}B\in R^{n\times(k-l)}$ have $\beta_{i}\neq0$ whereas the last
$k-l-l_{2}$ have $\beta_{i}=0.$ The matrix $B$ is obtained by generating
\[
B=\left(
\begin{array}
[c]{cc}%
B_{1} & 0\\
0 & B_{2}%
\end{array}
\right)
\]
where $B_{1}=(z_{1}/||z_{1}||\cdots z_{l_{2}}/||z_{l_{2}}||)$ with
$z_{1},\ldots,z_{l_{2}}\overset{i.i.d.}{\sim}N_{l_{1}}(0,I)$ independent of
$B_{2}=(z_{l_{2}+1}/||z_{l_{2}+1}||\cdots z_{k-l-l_{2}}/||z_{l_{k-l-l_{2}}%
}||)$ with $z_{l_{2}+1},\ldots,z_{k-l-l_{2}}$ i.i.d. $N_{l-l}(0,I).$ Note that
this ensures that the columns of $X_{2}$ are all standardized. Furthermore,
since it is assumed that the last $l-l_{1}$ variables of $X_{1}$ and the last
$k-l-l_{2}$ variables of $X_{2}$ don't have an effect, the form of $B$ is
necessarily of the diagonal form given. For, if it was allowed that the last
$k-l-l_{2}$ columns of $X_{2}$ were linearly dependent on the the first
$l_{1}$ columns of $X_{1},$ then this would induce a dependence on the
corresponding variables and this is not the intention in the simulation.
Similarly, if the first $l_{2}$ columns of $X_{2}$ were dependent on the last
$l-l_{1}$ columns of $X_{1},$ then this would imply that the variables
associated with these columns of $X_{1}$ have an effect and this is not the intention.

The sampling model is then prescribed by setting $l=10,l_{1}=5,l_{2}=2,$ with
$\beta_{i}=4$ for $i=1,\ldots,5,11,12$ with the remaining $\beta_{i}=0,$
$\sigma^{2}=1,m=2,$ so $n=22$ and we consider various values of $k\geq l.$ It
is is to be noted that a different data set was generated for each value of
$k.$ The prior is specified as in Example 4 where the values $\lambda_{0}%
^{2}=4,\alpha_{1}=11,\alpha_{2}=12$ were chosen so that there will be no
prior-data conflict arising with the generated data. Also, we considered
several values for the discretization parameter $\delta.$ A hypothesis was
classified as true if the relative belief ratio is greater than 1 and
classified as false if it is less than 1. Table \ref{confusereg} gives the
confusion matrices with $\delta=0.1.$ The value $\delta=0.5$ was also
considered but there was no change in the results.

One fact stands out immediately, namely, in all of these example only one
misclassification was made and this was in the full rank ($k=10$) case where
one hypothesis which was true was classified as a positive. The effect sizes
that exist are reasonably large, and so it can't be expected that the same
performance will arise with much smaller effect sizes, but it is clear that
the approach is robust to the number of hypotheses considered. It should also
be noted, however, that the amount of data is relatively small and the success
of the procedure will only improve as this increases. This result can, in
part, be attributed to the fact that a logically sound measure of evidence is
being used. $\blacksquare$%

\begin{table}[tbp] \centering
\begin{tabular}
[c]{|c|}\hline%
\begin{tabular}
[c]{llll}%
$k=10$ & Classified Positive & Classified Negative & Total\\
True Positive & \multicolumn{1}{c}{5} & \multicolumn{1}{c}{0} &
\multicolumn{1}{c}{5}\\
True Negative & \multicolumn{1}{c}{1} & \multicolumn{1}{c}{4} &
\multicolumn{1}{c}{5}\\
Total & \multicolumn{1}{c}{6} & \multicolumn{1}{c}{4} & \multicolumn{1}{c}{10}%
\end{tabular}
\\\hline%
\begin{tabular}
[c]{llll}%
$k=20$ & Classified Positive & Classified Negative & Total\\
True Positive & \multicolumn{1}{c}{7} & \multicolumn{1}{c}{0} &
\multicolumn{1}{c}{7}\\
True Negative & \multicolumn{1}{c}{0} & \multicolumn{1}{c}{13} &
\multicolumn{1}{c}{13}\\
Total & \multicolumn{1}{c}{7} & \multicolumn{1}{c}{13} &
\multicolumn{1}{c}{20}%
\end{tabular}
\\\hline%
\begin{tabular}
[c]{llll}%
$k=50$ & Classified Positive & Classified Negative & Total\\
True Positive & \multicolumn{1}{c}{7} & \multicolumn{1}{c}{0} &
\multicolumn{1}{c}{7}\\
True Negative & \multicolumn{1}{c}{0} & \multicolumn{1}{c}{43} &
\multicolumn{1}{c}{43}\\
Total & \multicolumn{1}{c}{7} & \multicolumn{1}{c}{43} &
\multicolumn{1}{c}{50}%
\end{tabular}
\\\hline%
\begin{tabular}
[c]{llll}%
$k=100$ & Classified Positive & Classified Negative & Total\\
True Positive & \multicolumn{1}{c}{7} & \multicolumn{1}{c}{0} &
\multicolumn{1}{c}{7}\\
True Negative & \multicolumn{1}{c}{0} & \multicolumn{1}{c}{93} &
\multicolumn{1}{c}{93}\\
Total & \multicolumn{1}{c}{7} & \multicolumn{1}{c}{93} &
\multicolumn{1}{c}{100}%
\end{tabular}
\\\hline
\end{tabular}
\caption{Confusion matrices for the numerical example in Example 5.}\label{confusereg}%
\end{table}%

\section{Conclusions}

An approach to the problem of multiple testing has been developed based on the
relative belief ratio. It is argued in Evans (2015) that the relative belief
ratio is a valid measure of evidence as it measures change in belief as
opposed to belief and, among the many possible candidates for such a measure,
it is the simplest with the best properties. One can expect that statistical
procedures based on valid measures of evidence will perform better than
procedures that aren't as they possess a sounder logical basis. For the
multiple testing problem this is reflected in the increased flexibility as
well as in the results.

It seems that an appropriate measure of evidence in a statistical problem
requires the specification of a prior. While this may be controversial to
some, it is to be noted that there are tools for dealing with the subjective
nature of some of the ingredients to a statistical analysis such as the
sampling model and prior. In particular, there is the process of checking for
prior-data conflict after the data is obtained and possibly modifying the
prior based upon the idea of weak informativity when such a conflict is
encountered. Before the data is actually collected, one can measure to what
extent a particular prior will bias the results based upon the particular
measure of evidence used. If bias is encountered several mitigating steps can
be taken but primarily this will require increasing the amount of data
collected. These concepts play a key role in the multiple testing problem.

\section*{Acknowledgements}

Thanks to Professor Lei Sun for making her notes on multiple testing available.

\section*{Bibliography}

\noindent Baskurt, Z. and Evans, M. (2013) Hypothesis assessment and
inequalities for Bayes factors and relative belief ratios. Bayesian Analysis,
8, 3, 569-590.\smallskip

\noindent Cao, Y., Evans, M. and Guttman, I. (2014) Bayesian factor analysis
via concentration. Current Trends in Bayesian Methodology with Applications,
edited by S. K. Upadhyay, U. Singh, D. K. Dey and A. Loganathan. CRC
Press.\smallskip

\noindent Carvalho, C. M., Polson, N. G. and Scott, J. G. (2009) Handling
sparsity via the horseshoe. Journal of Machine Learning Research W\&CP 5:
73-80.\smallskip

\noindent Efron, B., Hastie, T, Johnstone, I., and Tibshirani, R. (2006) Least
angle regression. The Annals of Statistics, 32, 407-499.\smallskip

\noindent Evans, M. (1997) Bayesian inference procedures derived via the
concept of relative surprise. Communications in Statistics, 26,
1125-1143.\smallskip

\noindent Evans , M. (2015) Measuring Statistical Evidence Using Relative
Belief. Monographs on Statistics and Applied Probability 144, CRC
Press.\smallskip

\noindent Evans, M. and Jang, G. H. (2011a) A limit result for the prior
predictive. Statistics and Probability Letters, 81, 1034-1038.\smallskip

\noindent Evans, M. and Jang, G. H. (2011b) Weak informativity and the
information in one prior relative to another. Statistical Science, 26, 3,
423-439.\smallskip

\noindent Evans, M. and Moshonov, H. (2006) Checking for prior-data conflict.
Bayesian Analysis, 1, 4, 893-914.\smallskip

\noindent George, E. I. and McCulloch, R. E. (1993). Variable selection via
Gibbs sampling. Journal of the American Statistical\ Association,
88,881-889.\vspace{2pt}

\noindent George, E. I. and McCulloch, R. E. (1997). Approaches for Bayesian
variable selection. Statistica Sinica, 7(2), 339--373.\vspace{2pt}

\noindent Hastie, T., Tibshirani, R. and Martin Wainwright, M. (2015)
Statistical Learning with Sparsity: The Lasso and Generalizations. Monographs
on Statistics and Applied Probability 143, CRC Press.\smallskip

\noindent Park, R. and Casella, G. (2008) The Bayesian Lasso. Journal of the
American Statistical Association, 103, 681-686.\smallskip

\noindent Rockova, V. and George, E. I. EMVS: The EM approach to Bayesian
variable selection. Journal of the American Statistical Association, 109, 506,
828 - 846.\vspace{2pt}

\noindent Tibshirani, R. (1996). Regression shrinkage and selection via the
lasso. Journal of the Royal Statistical Society B., 58, 1, 267-288.

\section*{Appendix}

\noindent Proof of Lemma \ref{combin}: Let $\Delta(i)$ be the event that
exactly $i$ of $A_{1},\ldots,A_{k}\in\mathcal{F}$ occur, so that $\cup
_{i=1}^{k}A_{i}=\cup_{i=1}^{k}\Delta(i)$ and note that the $\Delta(i)$ are
mutually disjoint. When $l<k,$
\begin{align*}
S_{l,k}  &  =%
{\displaystyle\sum\limits_{\{i_{1},\ldots,i_{l}\}\subset\{1,\ldots,k\}}}
I_{A_{i_{1}}\cup\cdots\cup A_{i_{l}}}=\binom{k}{l}\sum_{i=0}^{l-1}%
I_{\Delta(k-i)}+\sum_{i=l}^{k-1}\left[  \binom{k}{l}-\binom{i}{l}\right]
I_{\Delta(k-i)}\\
&  =\binom{k}{l}\sum_{i=0}^{k-1}I_{\Delta(k-i)}-\sum_{i=l}^{k-1}\binom{i}%
{l}I_{\Delta(k-i)}%
\end{align*}
and $S_{k,k}=I_{A_{1}\cup\cdots\cup A_{k}}.$ Now consider $\binom{k}{l}%
^{-1}S_{l,k}-\binom{k}{l-1}^{-1}S_{l,k}$ which equals%
\begin{equation}
\frac{1}{\binom{k}{l}}%
{\displaystyle\sum\limits_{\{i_{1},\ldots,i_{l}\}\subset\{1,\ldots,k\}}}
I_{A_{i_{1}}\cup\cdots\cup A_{i_{l}}}-\frac{1}{\binom{k}{l-1}}%
{\displaystyle\sum\limits_{\{i_{1},\ldots,i_{l-1}\}\subset\{1,\ldots,k\}}}
I_{A_{i_{1}}\cup\cdots\cup A_{i_{l-1}}} \label{combin1}%
\end{equation}
If $l=k,$ then (\ref{combin1}) equals $I_{A_{1}\cup\cdots\cup A_{k}}%
-\sum_{i=0}^{k-1}I_{\Delta(k-i)}+I_{\Delta(1)}=I_{A_{1}\cup\cdots\cup A_{k}%
}-\sum_{i=0}^{k-2}I_{\Delta(k-i)}\ $which is nonnegative. If $l<k,$ then
(\ref{combin1}) equals $\binom{k}{l-1}^{-1}I_{\Delta(k-l+1)}+\sum_{i=l}%
^{k-1}[\binom{i}{l-1}\binom{k}{l-1}^{-1}-\binom{i}{l}\binom{k}{l}%
^{-1}]I_{\Delta(k-i)}$ which is nonnegative since an easy calculation gives
that each term in the second sum is nonnegative. The expectation of
(\ref{combin1}) is then nonnegative and this establishes the result.

\end{document}